\documentclass{elsarXiv}
\usepackage{natbib}
\usepackage{amssymb}
\usepackage{eucal}
\usepackage{amsmath}
\usepackage{amscd}
\usepackage{xspace}

\usepackage{multicol}

\usepackage{amsfonts,latexsym}
\usepackage[dvips,ps2pdf]{hyperref}
\usepackage{float}
\usepackage{fancybox}

\begin{document}  

\newcommand\notes[1]{{}}
\newcommand\deleted[1]{{}}
\newcommand\mkset[1]{{\{#1\}}}
\newcommand\paren[1]{{\rm(}#1{\rm)}}

\renewcommand\Bbb{\mathbb}  
\renewcommand\frak{\mathfrak} 

\newcommand{\binomial}[2]{\genfrac{(}{)}{0pt}{}{#1}{#2}}


\newcommand{\nc}{\newcommand}
\newcommand{\delete}[1]{}

\nc{\dfootnote}[1]{{ }}                   
\nc{\ffootnote}[1]{\dfootnote{#1}}

\delete{
\nc{\mfootnote}[1]{{}}        
\nc{\ofootnote}[1]{{}}        
}

\nc{\mfootnote}[1]{\footnote{#1}} 
\nc{\ofootnote}[1]{\footnote{\tiny Older version: #1}}



\nc{\mlabel}[1]{\label{#1}}  
\nc{\myear}[1]{(\citeyear{#1})}
\nc{\mauthor}[1]{\citeauthor{#1}}
\nc{\mcite}[1]{\citeauthor{#1}\,(\citeyear{#1})}

\nc{\mref}[1]{\ref{#1}}  
\nc{\mkeep}[1]{{}}      
\nc{\mbibitem}[1]{\bibitem{#1}} 

\delete{
\nc{\mcite}[1]{\cite{#1}{{\bf{{\ }(#1)}}}}
\nc{\mlabel}[1]{\label{#1}
{\hfill \hspace{1cm}{\bf{{\ }\hfill(#1)}}}}
\nc{\mref}[1]{\ref{#1}{{\bf{{\ }(#1)}}}}
\nc{\mbibitem}[1]{\bibitem[\bf #1]{#1}} 
\nc{\mkeep}[1]{\marginpar{{\bf #1}}} 
}

\deleted{
\newtheorem{theorem}{Theorem} [section]
\newtheorem{lemma}[theorem]{Lemma}
\newtheorem{coro}[theorem]{Corollary}
\newtheorem{prop-def}{Proposition-Definition}
\newtheorem{propprop}{Proposed Proposition}
\newtheorem{conjecture}{Conjecture}
\newtheorem{exam}{Example}[section]
\newtheorem{assumption}{Assumption}
\newtheorem{condition}[theorem]{Condition}
}

\newtheorem{theorem}{Theorem} [section]
\newtheorem{exam}[theorem]{Example}
\newtheorem{defi}[theorem]{Definition}

\newtheorem{coro}[theorem]{Corollary}
\newtheorem{lemma}[theorem]{Lemma}


\nc{\bond}{\vdash}

\nc{\dftimes}{\widetilde{\otimes}}
\nc{\dfl}{\succ}
\nc{\dfr}{\prec}
\nc{\dfc}{\circ}
\nc{\dfb}{\bullet}
\nc{\dft}{\star}
\nc{\dfcf}{{\mathbf k}}
\nc{\spr}{\cdot}
\nc{\disp}[1]{\displaystyle{#1}}
\nc{\bin}[2]{ (_{\stackrel{\scs{#1}}{\scs{#2}}})}  
\nc{\binc}[2]{ \left (\!\! \begin{array}{c} \scs{#1}\\
    \scs{#2} \end{array}\!\! \right )}  
\nc{\bbinc}[2]{ \left (\!\! \begin{array}{c} {#1}\\
    {#2} \end{array}\!\! \right )}  
\nc{\bincc}[2]{  \left ( {\scs{#1} \atop
    \vspace{-.5cm}\scs{#2}} \right )}  
\nc{\sarray}[2]{\begin{array}{c}#1 \vspace{.1cm}\\ \hline
    \vspace{-.35cm} \\ #2 \end{array}}
\nc{\bs}{\bar{S}}
\nc{\dcup}{\stackrel{\bullet}{\cup}}
\nc{\dbigcup}{\stackrel{\bullet}{\bigcup}}
\nc{\la}{\longrightarrow}
\nc{\fe}{\'{e}}
\nc{\rar}{\rightarrow}
\nc{\dar}{\downarrow}
\nc{\dap}[1]{\downarrow \rlap{$\scriptstyle{#1}$}}
\nc{\uap}[1]{\uparrow \rlap{$\scriptstyle{#1}$}}
\nc{\defeq}{\stackrel{\rm def}{=}}
\nc{\dis}[1]{\displaystyle{#1}}
\nc{\dotcup}{\ \displaystyle{\bigcup^\bullet}\ }
\nc{\hcm}{\ \hat{,}\ }
\nc{\hcirc}{\hat{\circ}}
\nc{\hts}{\hat{\shpr}}
\nc{\lts}{\stackrel{\leftarrow}{\shpr}}
\nc{\rts}{\stackrel{\rightarrow}{\shpr}}
\nc{\lleft}{[}
\nc{\lright}{]}
\nc{\uni}[1]{\tilde{#1}}
\nc{\free}[1]{\bar{#1}}
\nc{\den}[1]{\check{#1}}
\nc{\lrpa}{\wr}
\nc{\curlyl}{\left \{ \begin{array}{c} {} \\ {} \end{array}
    \right . \!\!\!\!\!\!\!}
\nc{\curlyr}{ \!\!\!\!\!\!\!
    \left . \begin{array}{c} {} \\ {} \end{array}
    \right \} }
\nc{\longmid}{\left | \begin{array}{c} {} \\ {} \end{array}
    \right . \!\!\!\!\!\!\!}
\nc{\ot}{\otimes}
\nc{\ora}[1]{\stackrel{#1}{\rar}}
\nc{\ola}[1]{\stackrel{#1}{\la}}
\nc{\scs}[1]{\scriptstyle{#1}}
\nc{\mrm}[1]{{\rm #1}}
\nc{\margin}[1]{\marginpar{\rm #1}}   
\nc{\dirlim}{\displaystyle{\lim_{\longrightarrow}}\,}
\nc{\invlim}{\displaystyle{\lim_{\longleftarrow}}\,}
\nc{\mvp}{\vspace{0.5cm}}
\nc{\svp}{\vspace{2cm}}
\nc{\vp}{\vspace{8cm}}
\nc{\proofbegin}{\noindent{\bf Proof: }}
\nc{\proofend}{$\square$ \vspace{0.5cm}}
\nc{\sha}{{\mbox{\cyr X}}}  
\nc{\ncsha}{{\mbox{\cyr X}^{\mathrm NC}}}
\nc{\ncshao}{{\mbox{\cyr X}^{\mathrm NC,\,0}}}
\nc{\shpr}{\diamond}    
\nc{\shpro}{\diamond^0}    
\nc{\shpru}{\check{\diamond}}
\nc{\catpr}{\diamond_l}
\nc{\rcatpr}{\diamond_r}
\nc{\lapr}{\diamond_a}
\nc{\lepr}{\diamond_e}
\nc{\vep}{\varepsilon}
\nc{\labs}{\mid\!}
\nc{\rabs}{\!\mid}
\nc{\hsha}{\widehat{\sha}}
\nc{\lsha}{\stackrel{\leftarrow}{\sha}}
\nc{\rsha}{\stackrel{\rightarrow}{\sha}}
\nc{\lc}{\lfloor}
\nc{\rc}{\rfloor}
\nc{\rbset}{R}
\nc{\rbnum}{r}
\nc{\rbfun}{\mathbf{R}}
\nc{\pset}{P}
\nc{\pnum}{p}
\nc{\pfun}{\mathbf{P}}
\nc{\spset}{SP}
\nc{\spnum}{sp}
\nc{\spgen}{\mathbf{SP}}

\nc{\ann}{\mrm{ann}}
\nc{\Aut}{\mrm{Aut}}
\nc{\can}{\mrm{can}}
\nc{\colim}{\mrm{colim}}
\nc{\Cont}{\mrm{Cont}}
\nc{\rchar}{\mrm{char}}
\nc{\cok}{\mrm{coker}}
\nc{\dtf}{{R-{\rm tf}}}
\nc{\dtor}{{R-{\rm tor}}}
\renewcommand{\det}{\mrm{det}}
\nc{\Div}{{\mrm Div}}
\nc{\End}{\mrm{End}}
\nc{\Ext}{\mrm{Ext}}
\nc{\Fil}{\mrm{Fil}}
\nc{\Frob}{\mrm{Frob}}
\nc{\Gal}{\mrm{Gal}}
\nc{\GL}{\mrm{GL}}
\nc{\Hom}{\mrm{Hom}}
\nc{\hsr}{\mrm{H}}
\nc{\hpol}{\mrm{HP}}
\nc{\id}{\mrm{id}}
\nc{\im}{\mrm{im}}
\nc{\incl}{\mrm{incl}}
\nc{\length}{\mrm{length}}
\nc{\LR}{\mrm{LR}}
\nc{\mchar}{\rm char}
\nc{\NC}{\mrm{NC}}
\nc{\mpart}{\mrm{part}}
\nc{\ql}{{\QQ_\ell}}
\nc{\qp}{{\QQ_p}}
\nc{\rank}{\mrm{rank}}
\nc{\rcot}{\mrm{cot}}
\nc{\rdef}{\mrm{def}}
\nc{\rdiv}{{\rm div}}
\nc{\rtf}{{\rm tf}}
\nc{\rtor}{{\rm tor}}
\nc{\res}{\mrm{res}}
\nc{\SL}{\mrm{SL}}
\nc{\Spec}{\mrm{Spec}}
\nc{\tor}{\mrm{tor}}
\nc{\Tr}{\mrm{Tr}}
\nc{\tr}{\mrm{tr}}

\nc{\ab}{\mathbf{Ab}}
\nc{\Alg}{\mathbf{Alg}}
\nc{\Algo}{\mathbf{Alg}^0}
\nc{\Bax}{\mathbf{Bax}}
\nc{\Baxo}{\mathbf{Bax}^0}
\nc{\RBo}{\mathbf{RB}^0}
\nc{\BRB}{\mathbf{RB}}
\nc{\Dend}{\mathbf{DD}}
\nc{\bfk}{{\bf k}}
\nc{\bfone}{{\bf 1}}
\nc{\base}[1]{{a_{#1}}}
\nc{\detail}{\marginpar{\bf More detail}
    \noindent{\bf Need more detail!}
    \svp}
\nc{\Diff}{\mathbf{Diff}}
\nc{\gap}{\marginpar{\bf Incomplete}\noindent{\bf Incomplete!!}
    \svp}
\nc{\FMod}{\mathbf{FMod}}
\nc{\RB}{\mathbf{RB}}
\nc{\Int}{\mathbf{Int}}
\nc{\Mon}{\mathbf{Mon}}
\nc{\remarks}{\noindent{\bf Remarks: }}
\nc{\Rep}{\mathbf{Rep}}
\nc{\Rings}{\mathbf{Rings}}
\nc{\Sets}{\mathbf{Sets}}
\nc{\DT}{\mathbf{DT}}

\nc{\BA}{{\Bbb A}}
\nc{\CC}{{\Bbb C}}
\nc{\DD}{{\Bbb D}}
\nc{\EE}{{\Bbb E}}
\nc{\FF}{{\Bbb F}}
\nc{\GG}{{\Bbb G}}
\nc{\HH}{{\Bbb H}}
\nc{\LL}{{\Bbb L}}
\nc{\NN}{{\Bbb N}}
\nc{\QQ}{{\Bbb Q}}
\nc{\RR}{{\Bbb R}}
\nc{\TT}{{\Bbb T}}
\nc{\VV}{{\Bbb V}}
\nc{\ZZ}{{\Bbb Z}}


\nc{\cala}{{\mathcal A}}
\nc{\calc}{{\mathcal C}}
\nc{\cald}{{\mathcal D}}
\nc{\cale}{{\mathcal E}}
\nc{\calf}{{\mathcal F}}
\nc{\calg}{{\mathcal G}}
\nc{\calh}{{\mathcal H}}
\nc{\cali}{{\mathcal I}}
\nc{\calj}{{\mathcal J}}
\nc{\call}{{\mathcal L}}
\nc{\calm}{{\mathcal M}}
\nc{\caln}{{\mathcal N}}
\nc{\calo}{{\mathcal O}}
\nc{\calp}{{\mathcal P}}
\nc{\calr}{{\mathcal R}}
\nc{\calt}{{\mathcal T}}
\nc{\calw}{{\mathcal W}}
\nc{\calx}{{\mathcal X}}
\nc{\CA}{\mathcal{A}}

\nc{\fraka}{{\mathfrak a}}
\nc{\frakB}{{\mathfrak B}}
\nc{\frakb}{{\mathfrak b}}
\nc{\frakd}{{\mathfrak d}}
\nc{\frakF}{{\mathfrak F}}
\nc{\frakf}{{\mathfrak f}}
\nc{\frakg}{{\mathfrak g}}
\nc{\frakm}{{\mathfrak m}}
\nc{\frakM}{{\mathfrak M}}
\nc{\frakMo}{{\mathfrak M}^0}
\nc{\frakMl}{{\mathfrak M}^1}
\nc{\frakp}{{\mathfrak p}}
\nc{\frakx}{{\mathcal X}}
\nc{\ox}{\overline{\frakx}}
\nc{\frakX}{{\mathfrak X}}
\nc{\fraky}{{\mathfrak y}}

\renewcommand\geq{\geqslant}
\renewcommand\leq{\leqslant}
\renewcommand\bar[1]{\overline{#1}}

\nc\rbop{{\lc\,\,\rc}}
\nc\rbopi[1]{{\lc_{#1} \, \rc_{#1}}}
\newcommand\rbo[1]{{\lc\,\,\rc^{(#1)}}}

\nc\rseq[1]{r_{#1}}
\nc\aseq[1]{a_{#1}}
\nc\bseq[1]{b_{#1}}
\nc\dseq[1]{d_{#1}}
\nc\iseq[1]{i_{#1}}
\nc\rdseq[2]{r_{#1,#2}}
\nc\adseq[2]{a_{#1,#2}}
\nc\bdseq[2]{b_{#1,#2}}
\nc\ddseq[2]{d_{#1,#2}}
\nc\idseq[2]{i_{#1,#2}}

\nc{\exponent}{{exponent}\xspace}
\nc{\exponents}{{exponents}\xspace}
\nc{\pexp}{{u}\xspace}
\nc{\xexp}{{v}\xspace}
\nc{\pexpv}{\vec{u}\xspace}
\nc{\xexpv}{\vec{v}\xspace}

\nc{\xarity}{{m}\xspace}
\nc{\pdeg}{{n}\xspace}
\nc{\xarityv}{\vec{m}\xspace}
\nc{\pdegv}{\vec{n}\xspace}

\nc{\run}{run\xspace}
\nc{\runs}{runs\xspace}
\nc{\prun}{{k}\xspace}
\nc{\xrun}{{\ell}\xspace}
\nc{\prunv}{\vec{k}\xspace}
\nc{\pxrun}{\vec{\ell}\xspace}
\nc{\xb}{b\xspace}
\nc{\pa}{{a}\xspace}

\nc{\prunlen}{{\mu}\xspace}
\nc{\xrunlen}{{\nu}\xspace}

\nc{\pindex}{{\alpha}\xspace}
\nc{\xindex}{{\beta}\xspace}

\nc{\pdegvar}{{z}\xspace}
\nc{\xarityvar}{{t}\xspace}
\nc{\prunvar}{{\zeta}\xspace}
\nc{\xrunvar}{{\theta}\xspace}
\nc{\pdegvarv}{\vec{z}\xspace}
\nc{\xarityvarv}{\vec{t}\xspace}
\nc{\prunvarv}{\vec{\zeta}\xspace}
\nc{\xrunvarv}{\vec{\theta}\xspace}
\nc{\inftyv}{\stackrel{\rightarrow}{\infty}\xspace}
\nc{\onev}{\stackrel{\rightarrow}{1}\xspace}
\nc{\comGF}{{\mathbf G}\xspace}
\nc{\comfun}{{g}\xspace}
\nc{\comset}{{G}\xspace}
\nc{\xnumb}{{q}\xspace}
\nc{\pnumb}{{p}\xspace}

\nc{\pv}{{\vec{P}}}
\nc{\xv}{{\vec{x}}}
\nc{\oX}{\overline{X}}

\nc{\pxmap}{{\Phi}}
\nc{\pmap}{{\Phi^\pv}}
\nc{\xmap}{{\Phi^\xv}}
\nc{\pxmappart}{{\Phi^{\pv', \xv'}}}
\nc{\pxmaplong}{{\Phi_{\pexpv, \xexpv, \pdegv, \xarityv}}}
\nc{\pxmapscalar}{{\Phi_{\pexpv, \xexpv, \pdeg, \xarity}}}
\nc{\pmaplong}{{\Phi_{\pexpv, \xexpv, \pdegv, \xarityv}^\pv}}
\nc{\xmaplong}{{\Phi_{\pexpv, \xexpv, \pdegv, \xarityv}^\xv}}
\nc{\pxmappartlong}{{\Phi_{\pexpv, \xexpv, \pdegv, \xarityv}^{\pv',\xv'}}}

\nc{\pxrsm}{{\Psi}}
\nc{\prsm}{{\Psi^\pv}}
\nc{\xrsm}{{\Psi^\xv}}
\nc{\pxrsmpart}{{\Psi^{\pv', \xv'}}}
\nc{\pxrsmlong}{{\Psi_{\pexpv, \xexpv, \pdegv, \xarityv}}}
\nc{\prsmlong}{{\Psi_{\pexpv, \xexpv, \pdegv, \xarityv}^\pv}}
\nc{\xrsmlong}{{\Psi_{\pexpv, \xexpv, \pdegv, \xarityv}^\xv}}
\nc{\pxrsmpartlong}{{\Psi_{\pexpv, \xexpv, \pdegv, \xarityv}^{\pv',\xv'}}}
\nc{\pxscalar}{{\Psi_{\pexp, \xexp, \pdeg, \xarity}}}
\nc{\pxinfty}{{\Psi_{\infty, \infty, \pdeg, \xarity}}}
\nc{\pscalar}{{\Psi_{\pexp, \xexp, \pdeg, \xarity}^P}}
\nc{\xscalar}{{\Psi_{\pexp, \xexp, \pdeg, \xarity}^x}}

\nc{\pxvec}{{\Psi_{\pexpv, \xexpv, \pdeg, \xarity}}}

\nc\colorset{{C}}
\nc\colornum{{c}}
\nc\colorGF{{\mathbf C}}
\nc\F{{\mathbf F}}
\font\cyr=wncyr10

\nc{\redtext}[1]{\textcolor{red}{#1}}

\newcommand\RBA{{Rota-Baxter algebra}\xspace}
\newcommand\RBAs{{Rota-Baxter algebras}\xspace}
\newcommand\round{{\central}\xspace}
\newcommand\paired{{indecomposable}\xspace}
\newcommand\unpaired{{decomposable}\xspace}
\newcommand\central{{bracketed}\xspace}
\newcommand\associates{{associates}\xspace}
\newcommand\associate{{associate}\xspace}
\newcommand\RBW{{\tt RBW}\xspace}
\newcommand\rbw{{\tt RBW}\xspace}
\newcommand\RBWs{{\tt RBWs}\xspace}
\newcommand\rbws{{\tt RBWs}\xspace}
\newcommand\angled[1]{{\langle {\tt #1} \rangle}\xspace}
\newcommand\bround{\angled{\round}\xspace}
\newcommand\bRBW{\angled{\RBW}\xspace}
\newcommand\bpaired{\angled{\paired}\xspace}
\newcommand\bunpaired{\angled{\unpaired}\xspace}
\newcommand\bassociate{\angled{\associate}\xspace}
\newcommand\bempty{{\emptyset}\xspace}

\newfont{\xn}{cmr10 scaled 800}
\newfont{\xm}{cmr10 scaled 640}

\begin{frontmatter}



\title{Enumeration of Rota-Baxter Words{\raise1ex\hbox{\small\dag}}}

\author[li]{Li Guo} and
\author[sit]{William Y.
Sit\corauthref{cor1}}
\address[li]{Dept. of Math. \& Comp. Sci.,
         Rutgers University,
         Newark, NJ 07102}
\ead{liguo@newark.rutgers.edu (Li Guo)}
\ead[url]{http://newark.rutgers.edu/\~{ }liguo (Li Guo)}
\address[sit]{Dept. of Math., City College, City
University of New York, New York, NY 10031}
\ead{wyscc@sci.ccny.cuny.edu}
\ead[url]{http://scisun.sci.ccny.cuny.edu/\~{ }wyscc}

\corauth[cor1]{Corresponding author.
}

\begin{abstract}
\setlength{\baselineskip}{12pt}
\setlength{\normalbaselineskip}{12pt}
In this paper, we prove results on enumerations of sets of Rota-Baxter
words in a finite number of generators and a finite number of unary
operators.  Rota-Baxter words are words formed by concatenating
generators and images of words under Rota-Baxter operators.  Under
suitable conditions, they form canonical bases of free Rota-Baxter
algebras and are studied recently in relation to combinatorics, number
theory, renormalization in quantum field theory, and operads.
Enumeration of a basis is often a first step to choosing a data
representation in implementation.  Our method applies some simple
ideas from formal languages and compositions (ordered partitions) of
an integer.  We first settle the case of one generator and one
operator where both have exponent 1 (the idempotent case).  Some
integer sequences related to these sets of Rota-Baxter words are known
and connected to other combinatorial sequences, such as the Catalan
numbers, and others are new.  The recurrences satisfied by the
generating series of these sequences prompt us to discover an
efficient algorithm to enumerate the canonical basis of certain free
Rota-Baxter algebras.  More general sets of Rota-Baxter words are
enumerated with summation techniques related to compositions of
integers.
\end{abstract}

\begin{keyword}Enumerative Combinatorics \sep
Rota-Baxter words and algebras\sep grammar \sep Catalan
numbers \sep generating functions \sep compositions.

{\it MSC codes:} 16W99, 05A15, 08B20, 68W30.
\notes{ 16W99 (algebra with extra structures), 05A15(enumeration and
generating functions), 08B20(free algebras in universal algebras),
68W30 (computational algebra)}


\end{keyword}
\end{frontmatter}
\vfill
$\,^\dag${This
revised full version of the extended abstract published in
Proc. ISSAC 2006 corrected minor errors and strengthened some
results.}


\setcounter{section}{0}

\vfill \section{Introduction}

This paper studies enumeration and algorithms related to generation of
sets of Rota-Baxter words which occur naturally as canonical bases of
certain Rota-Baxter algebras.

In the 1950s, Spitzer proved a fundamental identity on fluctuation
theory in probability by analytic methods.  The field of Rota-Baxter
algebra was started after G.~Baxter~\myear{Ba} showed that Spitzer's
identity follows more generally by a purely algebraic argument for any
linear operator $P$ on an algebra satisfying the simple identity
\begin{equation}
P(x)P(y)=P(xP(y)+P(x)y + \lambda xy)
\mlabel{eq:RB}
\end{equation}
for all elements $x, y$ of the algebra, where $\lambda$ is a
constant.\footnote{We will take $\lambda=-1$ in this paper.  Some
authors use $-\theta$ in place of $\lambda$.} Rota studied this
operator through his many articles and communications
(see~\mcite{Ro2}, for example).  As a recognition to the great
contribution of Rota and Baxter, such a linear operator $P$ is called
a {\bf Rota-Baxter operator} and an algebra with the operator is
called a {\bf Rota-Baxter algebra}.  In spite of diverse applications
of \RBAs in mathematics and physics, the study of \RBA itself has been
highly combinatorial.  \mcite{Ro1} and \mcite{R-S}, for instance,
related Rota-Baxter operator to other combinatorial identities, such
as the Waring formula and Bohnenblust-Spitzer identity.  Explicit
constructions of free commutative \RBAs have played an important role
in further studies, from Cartier~\myear{Ca} and Rota~\myear{Ro1} in
the 1970s, to \mauthor{G-K} \myear{G-K,G-K2} in the 1990s.  Because of
the combinatorial nature of the constructions, the related enumeration
problems are interesting to study.  For example, \mcite{Gu} showed
that free commutative \RBAs on the empty set are related to Stirling
numbers of the first and second kind, and these results in general
provide generating series for other number sequences.

The unexpected application of non-commutative \RBAs in the work of
Connes and Kreimer\!  \myear{C-K1,C-K2} and Ebrahimi-Fard,\!  Guo and
Kreimer \myear{E-G-K3,E-G-K2} on renormalization of quantum field
theory moves the constructions of the corresponding free objects to
the forefront.  Such constructions were obtained recently
by~\mauthor{E-G}\,\myear{E-G,E-G2}, providing a fuller understanding
of the connection first made by Aguiar (see \mcite{Ag} and \mcite{EF})
between Rota-Baxter algebras and dendriform algebras of Loday and
Ronco~\myear{Lo}, and in particular their Hopf algebra of planar
trees.  Free nonunitary Rota-Baxter algebras are convenient in the
study of the adjoint functor from dendriform algebras to Rota-Baxter
algebras.

We consider enumeration and algorithmic generation of sets of strings
called Rota Baxter words (\rbws) which represent expressions in
generators and unary operators.  Under suitable conditions, these sets
form the canonical bases for (non-unitary) free Rota-Baxter algebras
with finitely many operators and finitely many generators.  These
constructions of free Rota-Baxter algebras are only recently explored
in special cases in \mauthor{E-G} \myear{E-G, E-G2} and we enumerate
not only their canonical bases but also sets that may be useful for
more general yet-to-be-explored constructions.  For this purpose, we
apply concepts and methods from formal languages, grammars,
compositions of integers, and generating functions.  While most of the
sequences and double sequences are found to be closely related to the
Catalan numbers, we found sequences that are not covered in the Sloane
data base.  The enumeration study leads us to an algorithm that
generates canonical bases for certain free Rota-Baxter algebras.
These algorithmic explorations of enumeration methods (either
sequential or randomized) are the first steps that must precede
development of important software tools in symbolic computation
packages that allow further investigations of algebraic properties of
related algebras based on Rota-Baxter words.  In addition to
Rota-Baxter word sets in several free operators and free generators,
we consider sets of \rbws where the number of consecutive applications
of operators and the number of consecutive generators are bounded.
These \rbws contain interesting combinatorial structures related to
planar rooted trees which form Rota-Baxter algebras as developed
in~\mcite{E-G2}.  Some of these properties in certain cases are
studied by~\mcite{A-M}.  Specifically, Aguiar and
Moreira used bijections between combinatorial objects while we applied
direct algorithmic enumerations.

The rest of the paper is organized into sections dealing with three
levels of generalities with each level built upon the previous one.
Section \mref{sec:zero} briefly reviews concepts on a free Rota-Baxter
algebra and its canonical basis consisting of Rota-Baxter words.
Section \mref{sec:one} deals with the case of one idempotent operator
and one idempotent generator (that is, \exponent 1 case).  In Section
\mref{sec:two},
we generalize this to arbitrary \exponents for one operator and one
generator.  In Section \mref{sec:three}, we further generalize this to
arbitrary number of operators and generators, with arbitrary
uniform \exponents relations.  We end with
a brief remark on future research.

\section{Background and notations} \mlabel{sec:zero}

Let $\bfk$ be a commutative unitary ring and let $A$ be a non-unitary
$\bfk$-algebra.  A non-unitary {\bf free Rota-Baxter algebra over $A$}
is a Rota-Baxter algebra $(F(A),P_A)$ together with a $\bfk$-algebra
homomorphism $j_A:  A\to F(A)$ with the property that, for any
Rota-Baxter algebra $(R,P)$ together with a $\bfk$-algebra
homomorphism $f:  A\to R$, there is a unique homomorphism $\free{f}:
F(A)\to R$ of Rota-Baxter algebras such that $f=\free{f} \circ j_A$.

Ebrahimi-Fard and Guo\,\!\myear{E-G, E-G2} explicitly constructed a
free,\! non-commut-ative, Rota-Baxter algebra over $A$, denoted by
$\ncshao(A)$, in the case a $\bfk$-basis $X$ of the $\bfk$-algebra $A$
exists and is given.  Since the enumeration of a $\bfk$-basis of
$\ncshao(A)$ is the main subject of study of this paper, we briefly
recall its construction.  In what follows, the product of $x_1, x_2
\in X$ in the algebra $A$ is denoted by $x_1x_2$ or by $x_1 \cdot x_2$
if clarity is needed, and the repeated $n$-fold product of $x \in X$
in the algebra $A$ is denoted as usual by $x^n$.

Let $\lc$ and $\rc$ be symbols, called brackets, and let $X'=X\cup
\{\lc,\rc\}$.  Let $S(X')$ be the free (non-commutative) semigroup
generated by $X'$.  We can view an element $u \in S(X')$ as a string
made up of symbols from $X'$.  The product of two elements $u, v \in
S(X')$ is, by an abuse of notation, also denoted by the concatenation
$uv$ whenever there is no confusion, and by explicitly using the
concatenation operator $\sqcup$ as in $u \sqcup v$ otherwise.  It
should be emphasized that the operator $\sqcup$ is {\it not} a symbol
in $X'$ and is used solely for the purpose to resolve ambiguity in
cases when $u = x_1 \in X$ and $v = x_2 \in X$.  We will adopt the
convention that the notation $uv = x_1 x_2$ always means the product
$x_1 \cdot x_2$ in the algebra $A$ (and it may happen that $X$ is
closed under algebra multiplication so that $x_1x_2 \in X$), and the
concatenation of $x_1$ with $x_2$ as elements of $S(X')$ will always
be denoted by $x_1 \sqcup x_2$.  As we shall see, concatenation of two
(or more) elements of $X$ are explicitly excluded in the \RBW sets, in
particular, in any canonical basis of $\ncshao(A)$. So such usage
is very limited.

\begin{defi} {\rm\
A {\bf Rota-Baxter word (\RBW)} $w$ is an element of $S(X')$ that
satisfies the following conditions.
\vspace{-0.15in}
\begin{enumerate}
\item The number of $\lc$ in $w$ equals the number of $\rc$
in $w$; \mlabel{it:equalao}
\item Counting from the left to the right, the cumulative number
of $\lc$ at each location is always greater or equal to that of
$\rc$; \mlabel{it:leqao}
\item No subword
$x_1 \sqcup x_2$ occurs in $w$,
for any $x_1,x_2\in X$; \mlabel{it:x}
\item No subword $\rc\lc$ or $\lc\,\rc$ occurs in $w$;
\mlabel{it:ceilao}
\end{enumerate}
 }
\mlabel{de:apwao}
\end{defi}
Interpreting Definition \mref{de:apwao}, a Rota-Baxter word $w$ can be
represented uniquely by a finite string composed of one or more
elements of $X$, separated (if more than one $x$) by a left bracket
$\lc$ or by a right bracket $\rc$, where the set of brackets formed
balanced pairs, but neither the string $\rc \lc $ nor the string $\lc
\,\rc $ appears as a substring.  For example, when $X = \mkset{x}$,
the word $w = \lc \lc x\rc x\lc x\rc \rc x\lc x\rc $ is an \RBW, but
$\lc x \sqcup x \rc$, $\lc x^2 \rc$, $\lc x\rc \lc x\rc $, $x\rc x\lc
x$, and $\lc x\lc \,\rc x\rc $ are not.  The number of balanced pairs
of brackets in an \RBW is called its {\bf degree}.  The degree of $w$
in the above example is 4.  \notes{ and its arity is 5.}

Let $\frakMo(X)$ be the set of Rota-Baxter words and let $\frakMl(X)$
be $\frakMo(X)$ with the empty (or trivial) word $\bempty$ adjoined.

\begin{exam}{\rm
Let $\bfk$ be a field. Let $A = \bfk[x]$ be the polynomial ring in one
indeterminate $x$ over $\bfk$.
Then $X = \mkset{x^n \mid n \in \NN}$ is a $\bfk$-basis. In this
case, if $a + b = n$, then $\lc x^n \rc = \lc x^a x^b \rc$ is
an RBW, but $\lc x^a \sqcup x^b \rc$ is not.}
\mlabel{ex:unipoly}
\end{exam}

\begin{exam}
{\rm
Let $B = \bfk[x]/{\frak a}$ be the quotient ring of $A$ of Example
\mref{ex:unipoly} by the ideal $\frak a$ generated by $x^2-x$. Then
writing $\bar{1} = 1 + {\frak a}$ and $\bar{x} = x + {\frak a}$, the
set $X
= \mkset{\bar{1},\bar{x}}$ is a $\bfk$-basis of $B$. Here $\lc \bar{1}
\rc$, $\lc
\bar{1} \cdot \bar{x} \rc$, and $\bar{1} \lc \bar{x} \rc$ are \rbws
but $\lc \bar{1} \sqcup \bar{x} \rc$ is not.
}
\mlabel{ex:idem}
\end{exam}

Let $\ncshao(A)$ be the free $\bfk$-module with basis $\frakMo(X)$.
Ebrahimi-Fard and Guo\,\myear{E-G,E-G2} show that
 the following properties
\begin{equation}
\begin{array}{rcl}
x\diamond x' &=& x \cdot x'\\
x \diamond \lc w \rc & = &x \lc w \rc\\
\lc w \rc \diamond x &=& \lc w \rc x \\
\lc w \rc \diamond \lc w' \rc &=&
\lc \lc w \rc \diamond w'\rc +\lc w \diamond \lc w'\rc \rc
+\lambda \lc w \diamond w' \rc
\end{array}
\mlabel{eq:shprod0}
\end{equation}
for all $x, x' \in X$ and all $w, w' \in
\frakMo(X)$ uniquely
define an associative bilinear product $\diamond$ on $\ncshao(A)$.
This product, together with the linear
operator \begin{equation}
P_A: \ncshao(A) \to \ncshao(A),\qquad
    P_A(w)=\lc w\rc \ {\rm if}\ w\in \frakMo(X),
    \mlabel{eq:rbop}
\end{equation}
and the natural embedding
$$j_A: A \to \ncshao(A),\qquad j_A(x) = x\ {\rm if}\  x\in X,$$
makes $\ncshao(A)$ the free (non-unitary) Rota-Baxter algebra over $A$.
We will not need to know the explicit construction of this associative
bilinear product
in $\ncshao(A)$ for the rest of the paper.  However, we note that as
an element of the algebra $\ncshao(A)$, the string $\lc w \rc$ may be
interpreted as the image of the operator $P_A$ on $w$ for any $w \in
\frakMo(X)$ and that for any such $w$ writable as the concatenation
$uv$ for $u, v \in \frakMo(X)$, the concatenation can be viewed
as $u \diamond v$ (the first three cases of Eq. (\mref{eq:shprod0}))
and this justifies the abuse of notation and convention in
using concatenation for both the algebra multipication in $A$ and the
semigroup product in $S(X')$.

\begin{exam}{\rm
Let $\bfk$ be a field.  Let $A = \bfk \langle x_1, \dots, x_q \rangle$
be the polynomial ring in $q$ non-commutating indeterminates $x_1,
\dots, x_q$ over $\bfk$.  Then the set $X$ of (non-commutative)
monomials is a $\bfk$-basis of $A$.  If $q \geq 2$, then $\lc x_2^3
x_1^4 x_2^2 \rc = \lc x_2^3 \cdot x_1^4 \cdot x_2^2 \rc$ is an RBW,
but $\lc x_2^3 \sqcup x_1^4 \sqcup x_2^2 \rc$ is not.  Moreover, in
$\ncshao(A)$, we have ${1} \diamond {x_i} = {x_i} = {x_i} \diamond
{1}$ and ${1} \diamond \lc w \rc = {1} \lc w \rc$ for any $w \in
S(X')$.  } \mlabel{ex:three} \end{exam}

\begin{defi}{\rm The free (non-unitary) Rota-Baxter algebra
$\ncshao(A)$ of
Example \mref{ex:three} will be denoted by $\ncshao(q)$ and referred
to as the free (non-unitary) Rota-Baxter algebra on $q$ generators
$x_1, \dots, x_q$.  The corresponding $\bfk$-basis
$\frakMo(A)$ consisting of \rbws built from $X$ will be denoted by
$\frakMo(q)$. Any non-commutative monomial $x \in X$ will be called an
$x$-{\bf \run}. For any \rbw $w \in \frakMo(q)$, the {\bf arity}
of $w$ is the number of $x_1, \dots, x_q$ appearing in $w$, counted
with multiplicities. For example, the arity of an $x$-\run is the
total degree of the monomial it represents and the \rbw
$w = x_1x_2^2\lc x_2^3 x_1^4 x_2^2 \rc$ has two $x$-\runs and arity
12.  } \mlabel{de:rbgen} \end{defi}

\begin{exam}
{\rm
Let $A = \bfk \langle x_1, \dots, x_q \rangle$ be as in Example
\mref{ex:three},
let $\xexpv = (\xexp_1, \dots, \xexp_q)$ be a vector of $q$ positive
integers,
and $\frak a$ be the bilateral ideal of $A$ generated by the
polynomials
$x_i^{\xexp_i + 1} - x_i$, $1 \leq i \leq q$.  Let $B$ be the quotient
$\bfk$-algebra $A/{\frak a}$.
Writing $\bar{1} = 1 + {\frak a}$ and $\bar{x}_i = x_i + {\frak a}$,
let $X$ be the set
consisting of all non-commutative finite power products $$x =
\bar{x}_{j_1}^{e_{j_1}} \,
\cdots \,
\bar{x}_{j_{\ell-1}}^{e_{j_{\ell-1}}} \,
\bar{x}_{j_\ell}^{e_{j_\ell}} \,
\bar{x}_{j_{\ell+1}}^{e_{j_{\ell+1}}} \,
\cdots \,
\bar{x}_{j_{r}}^{e_{j_r}}
$$
  in
$\bar{x}_i$\,$(1\! \leq\! i \leq\! q)$, where the indices satisfy
$j_{\ell-1}
\ne j_{\ell}$ and $j_\ell \ne j_{\ell+1}$ for all\! $\ell$ $(2 \leq
\ell \leq
r-1)$, and the exponents satisfy $1 \leq e_{j_\ell} \leq
\xexp_{j_\ell}$
for $\ell$ $(1 \leq \ell \leq r)$. Then $X$ is a $\bfk$-basis of
$B$ and $\ncshao(B)$ is a free Rota-Baxter algebra on $B$.
} \mlabel{ex:four}
\end{exam}

\begin{defi}{\rm The free (non-unitary) Rota-Baxter algebra
$\ncshao(B)$ of Example \mref{ex:four} will be denoted by $\ncshao(q,
\xexpv)$ and be referred to as the free (non-unitary) Rota-Baxter
algebra on $q$ generators $\bar{x}_1, \dots, \bar{x}_q$ with \exponent
vector $\xexpv$.  The corresponding $\bfk$-basis $\frakMo(B)$
consisting of \rbws built from $X$ will be denoted by $\frakMo(q,
\xexpv)$.  Any non-commutative monomial $x \in X$ will be called an
$x$-{\bf \run}.  For any \rbw $w \in \frakMo(q, \xexpv)$, the {\bf
arity} of $w$ is the number of $\bar{x}_1, \dots, \bar{x}_q$ appearing
in a canonical representation of $w$ as an element of $B$, counted
with multiplicities.  For example, the arity of an $x$-\run is the
total degree of the monomial it represents and the \rbw $w =
\bar{x}_1\bar{x}_2^2\lc \bar{x}_2^3 \bar{x}_1^4 \bar{x}_2^2 \rc$ has
two $x$-\runs and arity 12, provided $\xexp_1 \geq 4$ and $\xexp_2
\geq 3$.} \mlabel{de:rbexp} \end{defi}

In this paper, we will enumerate \rbws in $\frakMo(q, \xexpv)$
(actually $\frakMl(q, \xexpv)$, after adjoining the empty \rbw
$\bempty$) with a given degree and arity, by giving algorithms to
generate them and generating functions that count them.  We begin with
$q=1$ and $v_1 = 1$ in Section \mref{sec:one} under some extra
hypothesis by restricting to a subset of \rbws, but generalize the
results to arbitrary $q$ and \rbws involving multiple unary operators.
For some of these generalizations, we note that the corresponding free
Rota-Baxter algebras have not been constructed and the enumeration of
the sets of \rbws is included for possible future applications.

\section{One idempotent operator and one idempotent generator
case} \mlabel{sec:one}

In this section, we restrict ourselves to Example \mref{ex:four} when
$q = 1$ and $\xexpv = (1)$, that is, $\bar{x}$ is idempotent, and we
further assume that the Rota-Baxter operator $P$ is also idempotent
(that is, $P(P(w)) = P(w)$ for all $w$).  These restrictions allow us
to first focus on the word structures of free Rota-Baxter algebra
constructions before considering other factors involved in more
general Rota-Baxter words.  Interestingly, in most applications of
Rota-Baxter algebra in quantum field theory, the operators are
idempotent.

We contribute three results for the enumeration of a
canonical basis of the free Rota-Baxter algebra in this special case.
After reviewing some preliminary material and setting up notations, we
consider generating functions based on the degree of the Rota-Baxter
words in Section \mref{sec:onepttwo}.  In Section \mref{sec:oneptthree}, we
refine the study to consider generating functions based on the degree
and arity.  In Sections \mref{sec:oneptfour}, we give an algorithm to
generate this canonical basis with given degree and arity.

For simplicity, we will drop the bar notation above the generator $x$.
Under our current hypothesis that both the single generator $x$ and
the operator $P$ are idempotent, let $\rbset=\rbset_{1,1}$ be the
subset of $\frakMl(X)$ consisting of $\bempty$ and Rota-Baxter words
$w$ composed of $x$'s and pairs of balanced brackets such that no two
$x$'s are adjacent, and no two pairs of brackets can be immediately
adjacent or nested.  In other words, the strings $\rc \lc $, $\lc
\,\rc $, and strings of the form $\lc \lc * \rc \rc $ where the
brackets are balanced pairs and where $*$ may be any \RBW, do not
appear as substrings of $w$.  For example, the \rbw $\lc \lc x\rc \rc
x$ is not element of $\rbset$.  For the rest of this section, all
\rbws are assumed to be in $\rbset$.

\subsection{Generating functions of one variable}
\mlabel{sec:onepttwo}

Let $R(\pdeg)$ be the subset of $\rbset$ of degree $\pdeg$ (with our
convention, $R(0) = \mkset{\bempty,x}$).  For $n > 0$, let $B(\pdeg)$
be the subset of $R(\pdeg)$ consisting of \RBWs that begin with a left
bracket and end with a right bracket.  Words in $B(n)$ are said to be
{\bf \round}.
By pre- or post- concatenating a \round \RBW $w$ with $x$, we get
three new \RBWs:  $x w$, $w x$, and $x w x$, which are called
respectively the {\bf left}, {\bf right}, and {\bf bilateral
\associate} of $w$.  We also consider $x$ to be an associate of
the trivial word $\bempty$. Any non-trivial \RBW is either
\round or an \associate.  Thus for $n > 0$, the set $A(n)$ of all
\associates form the complement of $B(n)$ in $R(n)$ and it is the
disjoint union of these cosets:
\vspace{-0.05in}
$$A(n) = xB(n) \cup
B(n) x \cup xB(n) x, \qquad n > 0.
\vspace{-0.05in}$$
The set of \round \RBWs is further divided into two disjoint subsets.
The first subset $I(n)$ consists of all {\bf \paired} \round \RBWs,
whose beginning left bracket and ending right bracket are paired.  The
second subset $D(n)$ consists of all {\bf \unpaired} \round \RBWs
whose beginning left bracket and ending right bracket are not paired.
For convenience in counting, we define $B(0)$, $I(0)$, $D(0)$ to
be the empty set and note that
$A(0)\!  =\!  \mkset{x}$.

The following table lists these various types of \rbws in lower degrees.
\smallskip
\begin{center}
\begin{tabular}{| c | c | c || c | c  |} \hline
$\deg$ & $I(n)$ & $D(n)$ & $A(n)$ & $B(n)$ \\
\hline \hline
0 &  & & x & \\
\hline
1 & $\lc x\rc$ & & $x\lc x\rc, \lc x\rc x, x\lc x \rc x$
& $\lc x\rc$
 \\
\hline
2 & $\lc x\lc x\rc\rc, \lc \lc x\rc x\rc, \lc x\lc x\rc x\rc$
    & $\lc x\rc x\lc x\rc$ & 12 associates
& $I(2) \cup D(2)$\\
    \hline
\end{tabular}
\end{center}
\smallskip

In terms of formal languages, we start with an alphabet $\Sigma$
of {\bf terminals} consisting of a special symbol $\bempty$ and the
three symbols $\lc , x$, and $\rc$, a set of {\bf non-terminals}
consisting
of $\bround$, $\bpaired$, $\bunpaired$, $\bassociate$ and the sentence
symbol $\bRBW$. Let the production rules be:
\vspace{-0.2in}
\begin{eqnarray}
\bRBW &\rightarrow&  \bempty  \mid \bround \mid \bassociate
\mlabel{prod:rbw}\\
\bassociate &\rightarrow& x \mid x \bround \mid \bround x \mid x \bround x
\hspace{0.25in}
\mlabel{prod:assoc}\\
\bround &\rightarrow&  \bpaired \mid \bunpaired \mlabel{prod:round}\\
\hspace{-0.8in}\bpaired &\rightarrow& \lc \, \bunpaired \,\rc  \,\mid\, \lc \,
\bassociate \,\rc \mlabel{prod:paired}\\
\bunpaired &\rightarrow& \bround x \bround \mlabel{prod:unpaired}
\vspace{-0.2in}
\end{eqnarray}
By (\mref{prod:assoc}) and (\mref{prod:paired}), it is clear that the
sentences in this language will be \RBWs and vice versa. The
production rules thus define a {\bf grammar} whose {\bf language} will
be $\bigcup_{n=0}^\infty R(n)$. For the benefit of readers unfamiliar
with this area of computer science, we
refer them to \mcite{A-U} for basics and now give a detailed
proof.

Proof. Let ${\call}$ be the language defined by
the grammar above.  Clearly, $R(0) \subset \call$. It is easy to see
by induction on the length of
the string representing an \RBW $w$ that $w \in {\call}$ in case $w$
is either \paired or an associate since by Definition
\mref{de:apwao}, removing the outermost pair of balanced brackets from
an \paired \RBW or the appended $x$ (or $x$'s) from an
associate will yield an \RBW of shorter length.
If $w$ is \round and \unpaired, then we may write
$w = \lc w' \rc$ where the beginning and ending brackets are not
paired and $w' \in S(X')$.  Let $u$ be left subword of $w$ ending in
and including
the $\rc$ that balances the beginning $\lc$ of $w$. Then $u$ is
an \paired \RBW by definitions. Again by Definition
\mref{de:apwao}, the next symbol in $w$ following $u$ must be $x$,
which is then followed by a $\lc$ matching the ending $\rc$. Thus we
may write $w = u x v$ where $u, v$ are both \paired (and
\round) \RBWs of
shorter lengths and by induction, $u, v$, and hence also $w$ (by Rule
(\mref{prod:unpaired})), are in $\call$.

Conversely, to show that ${\call} \subseteq R$, it is only necessary
to show that a sentence production $\pi$ beginning with the rule
$\bRBW \rightarrow \bround$ will result in a string $w$ that is a
bracketed \RBW (by Rules (\mref{prod:rbw}), (\mref{prod:assoc}), and
Definition
\mref{de:apwao}). The shortest sentence production $\pi$ with this
property is $$\bRBW \rightarrow \bround \rightarrow \bpaired
\rightarrow \lc \bassociate \rc \rightarrow \lc x \rc$$
which results in an \RBW.  In general, we proceed by induction on the
number of productions applied in $\pi$.  The possible continuations of
$\bRBW \rightarrow \bround$ are those leading to $\bunpaired$, $\lc
\bunpaired \rc$, or $\lc \bassociate \rc$, each of which in turns
leads back to $\bround$.  Since application of any of the Rules
(\mref{prod:assoc}), (\mref{prod:paired}), or (\mref{prod:unpaired})
will result in an \RBW if the non-terminals involved result in \RBWs,
the induction hypothesis and Definition \mref{de:apwao} show that
these continuations will result in a valid \RBW.  \proofend

Because $A(n)$ and $B(n)$ are disjoint, and $I(n)$ and $D(n)$ are
disjoint, the proof shows that given any \RBW $w$, there is a unique
sentence production $\pi$ ending in $w$.  Sentences derivable from
$\bround$, $\bpaired$, $\bunpaired$, and $\bassociate$ correspond
bijectively respectively to \RBWs that are \round, \paired, \unpaired,
and \associate.

For $n \geq 0$, let $\rseq{n}$ (resp.~$\aseq{n}$, resp.~$\bseq{n}$,
resp.~$\iseq{n}$,
resp.~$\dseq{n}$) be the number of all (resp.~\associate,
resp.~\round,
resp.~\paired, resp.~\unpaired) \RBWs with $n$ pairs of (balanced)
brackets. The first few values of $\rseq{n}$ for $n = 0, 1, 2, \dots, 5$
are $$2, 4, 16, 80, 448, 2688, \dots $$ which suggests that it is the
sequence {\tt A025225} from \mcite{EIS}
 whose
$n$-th term is given by $2^{n+1} C_n$ . Here $C_n = \frac{1}{n+1}
\binomial{2n}{n}$ is the $n$-th Catalan number. By (\mref{prod:rbw}) and
(\mref{prod:assoc}), we clearly
have $\rseq{n} = 4\bseq{n}$ for $n>0$, and hence it suffices to prove
that $\bseq{n} = 2^{n-1}
C_n$. It is known that
the sequence A003645 whose $n$-th term is $2^{n-1} C_n$, $(n \geq
1)$, has a generating series given by
\begin{equation}
\sum_{n=1}^\infty
2^{n-1}C_n z^n = \frac{1-4z-\sqrt{1-8z}}{8z}.
\mlabel{eq:bseries}
\end{equation}
We will prove that $\bseq{n}$ indeed has a generating series given
by Eq.~\mref{eq:bseries}, and more.

\begin{theorem} \mlabel{thm:zero}
The generating series for $\rseq{n}, \bseq{n}, \iseq{n}, \dseq{n}$ and
$\aseq{n}$ are given by:
\begin{equation}
\begin{array}{rcl}
\rbfun(z)&=&\displaystyle \sum_{n=0}^\infty \rseq{n}z^n =
\frac{1-\sqrt{1-8z}}{2z}\,, \vspace{0.1in}\\
{\mathbf B}(z)&=&\displaystyle \sum_{n=0}^\infty \bseq{n} z^n =
\frac{1-4z-\sqrt{1-8z}}{8z}\,, \vspace{0.1in}\\
{\mathbf I}(z) &=& \displaystyle \sum_{n=0}^\infty \iseq{n}z^n =
\frac{1-2z-\sqrt{1-8z}}{2(z+1)}\,,\vspace{0.1in}\\
{\mathbf D}(z) &=&\displaystyle \sum_{n=0}^\infty \dseq{n}z^n =
\frac{1 - 7z + 4z^2 +(3z-1)\sqrt{1-8z}}{8z(z+1)}\,, \vspace{0.1in}\\
{\mathbf A}(z) &=& \displaystyle \sum_{n=0}^\infty \aseq{n}z^n =
\frac{3-4z -3\sqrt{1-8z}}{8z}.
\end{array}
\mlabel{thm:rbw}
\end{equation}
\end{theorem}

\begin{coro} \mlabel{thm:first}
The number of Rota-Baxter words of degree $n$ in the canonical basis
of the free Rota-Baxter algebra with a single idempotent generator and
idempotent operator is given by $$\rseq{n} = 2^{n+1} C_n, \qquad
n=0, 1, 2, \dots$$
where $C_n=\frac{1}{n+1} \binomial{2n}{n}$ is the $n$-th
Catalan number.
\end{coro}

\proofbegin
We note the initial values of the sequences are given by:
\begin{equation*}
\begin{array}{c}
\rseq{0}=2, \quad \bseq{0}=0, \quad \iseq{0} = 0, \quad \dseq{0}=0,
\quad \aseq{0} = 1, \vspace{-0.1in}\\
\rseq{1} = 4, \quad \bseq{1} = 1, \quad  \iseq{1} = 1, \quad \dseq{1}
= 0, \quad \aseq{1} = 3.
\end{array}
\end{equation*}
Thus, it does not matter whether the generating series for ${\mathbf
B}(z)$, ${\mathbf I}(z)$, and ${\mathbf D}(z)$ start at $n=0$ or
$n=1$.
By production rule (\mref{prod:round}), for all $n > 0$, we have
$\bseq{n}=\iseq{n}+\dseq{n}$. By production rule (\mref{prod:paired}),
an \RBW $w$ is in $I(n)$ if and only if $w =
\lc  w'\rc $  for some $w' \in A(n)$ or $w' \in D(n)$. Thus $\iseq{n}=
3\bseq{n-1} +
\dseq{n-1}$ and it follows that  $ \iseq{n}=4\bseq{n-1}-\iseq{n-1}$ or
$$ \iseq{n}+\iseq{n-1} =4 \bseq{n-1}.$$
Multiply both sides by $z^n$ and sum over $n\geq 2$. The left
hand side gives
$$\sum_{n= 2}^\infty \iseq{n} z^n +\sum_{n= 2}^\infty \iseq{n-1} z^n
={\mathbf I}(z)- z + z\,{\mathbf I}(z).$$
For the right hand side, note that by production rules
(\mref{prod:paired}) and (\mref{prod:round}), every \round, \unpaired
\RBW is formed from \round, \paired \RBWs with $x$ inserted between
them. This decomposition is clearly unique and the degree of the
\unpaired \RBW is the sum of the component \paired ones. Hence we have
\begin{equation}
 \dseq{n}=\sum_{\stackrel{(n_1,  \dots, n_p;\, n),}{p > 1}}
\iseq{n_1}\cdots \iseq{n_p} \mlabel{eq:recur}
 \end{equation}
where the notation $(n_1, \dots, n_p;\,n)$ denotes all compositions
$n_1 + \cdots + n_p = n$
of $n$ into $p$ positive integers, and the sum is over all
compositions for all lengths $p > 1$. Therefore,
noting the special case $p=1$ corresponds to a single summand $i_n$,
we have
\begin{equation}
 \bseq{n} = \iseq{n}+\dseq{n} =\sum_{\stackrel{(n_1,  \dots, n_p;\,
n),}{p \geq 1 }} \iseq{n_1}\cdots \iseq{n_p}\,, \quad (n \geq 1)
\mlabel{eq:decomp}  \end{equation}
{}From this, we obtain
$$
\begin{array}{rcl}
{\mathbf B}(z)&=& \displaystyle  \sum_{n= 1}^\infty \bseq{n} z^n
\vspace{0.1in}\\
&=& \displaystyle \sum_{n=1}^\infty \bigg( \sum_{(n_1,\dots,
n_p;\,
n),\,p\geq 1}\iseq{n_1}\cdots \iseq{n_p}\bigg) z^n \vspace{0.1in}\\
&=& \displaystyle \sum_{p= 1}^\infty \bigg (\sum_{\prun= 1}^\infty
\iseq{\prun}z^\prun\bigg)^p, \end{array}
$$
\begin{equation}{\mathbf B}(z) = \frac{{\mathbf I}(z)
}{1- {\mathbf I}(z)}\,.
\mlabel{eq:round} \end{equation}
Since
$$
4 \sum_{n= 2}^\infty \bseq{n-1} z^n = 4z\, {\mathbf B}(z)
$$
we have
$$ {\mathbf I}(z)-z+z\,{\mathbf I}(z)=\frac{4z\,{\mathbf
I}(z)}{1-{\mathbf I}(z)}.$$
Solving for ${\mathbf I}(z)$ with the initial condition $\iseq{0} =
0$, we have $$ {\mathbf I}(z)=\frac{1-2z-\sqrt{1-8z}}{2(z+1)}.$$
This proves the third formula of Eq.~(\mref{thm:rbw}) and by
Eq.~(\mref{eq:round}), after rationalizing the denominator,
we obtain the second formula of Eq.~(\mref{thm:rbw}).
The other formulas are easily derived from the relations
${\mathbf D}(z)\! =\! {\mathbf B}(z) - {\mathbf I}(z)$,
${\mathbf A}(z) = 1+3{\mathbf B}(z)$, and ${\mathbf R}(z)\! = \!
1+{\mathbf B}(z) + {\mathbf A}(z)$.  The corollary comes from
the remarks before the theorem.
\proofend

Theorem~\mref{thm:first} shows that the number $\iseq{n}$ of \round
\paired \RBWs
of degree $n$ ($n \geq 1$) is the $n$-th term of the sequence A062992:
$$1,3,13,67, 381, 2307, \dots$$
and the number $\dseq{n}$ of \round \unpaired \RBWs of
degree $n$ ($n \geq 2$) is a new sequence A115194, which starts with
$$1,7,45,291,1917,12867,\dots.$$

\subsection{Generating functions of two variables}
\mlabel{sec:oneptthree}

In our computational experiment, we observed that
the set $B(n)$, when stratified by the number of $x$'s appearing in an
\RBW, possesses certain very nice properties that may give better
combinatorial understanding of how the canonical basis is constructed
recursively (see the algorithm in the next subsection).
To describe the stratification, for any \RBW $w$, recall (Definition
\mref{de:rbexp}) that $w$ has {\bf arity} $\xarity$
if the number of $x$'s appearing in the string representation of $w$
is exactly $\xarity$.
For any $\xarity \geq 0$, let $R(n,\xarity)$ be
the subset of $\rbset$ of degree $n$ and arity $\xarity$, and define similarly
the notations $A(n,\xarity)$, $B(n,\xarity)$, $I(n,\xarity)$, and
$D(n,\xarity)$. These are all finite sets.
Let their sizes be respectively denoted by $\rdseq{n}{\xarity}$,
$\adseq{n}{\xarity}$, $\bdseq{n}{\xarity}$,
$\idseq{n}{\xarity}$, and $\ddseq{n}{\xarity}$. For initial values, we
have
\vspace{-0.1in}
\begin{equation*}
\begin{array}{rclcrcl}
\rdseq{0}{0} &=& 1; &\quad& \adseq{0}{0}&=& \bdseq{0}{0} =
\idseq{0}{0}=\ddseq{0}{0}=0; \vspace{-0.1in}\\
\rdseq{0}{1} &=& \adseq{0}{1}=1; &\quad&
\bdseq{0}{1}&=&\idseq{0}{1}=\ddseq{0}{1}=0; \vspace{-0.1in}
\\ \rdseq{1}{1} &=& \bdseq{1}{1} = \idseq{1}{1} = 1; &\quad&
\adseq{1}{1} &=&  \ddseq{1}{1} = 0; \vspace{-0.1in}\\
\rdseq{1}{2} &=& \adseq{1}{2} = 2; &\quad&
\bdseq{1}{2}&=&\idseq{1}{2}= \ddseq{1}{2} = 0; \vspace{-0.1in}\\
\rdseq{1}{3} &=& \adseq{1}{3} = 1; &\quad&
\bdseq{1}{3}&=&\idseq{1}{3}=\ddseq{1}{3} = 0;
\end{array}
\end{equation*}
\vspace{-0.2in}
\begin{equation*}
\begin{array}{rclcl}
\rdseq{0}{\xarity} &=& \adseq{0}{\xarity} =
\bdseq{0}{\xarity}=\idseq{0}{\xarity}=\ddseq{0}{\xarity}=0 &\quad&
{\rm\ for\ }\xarity \geq 2; \vspace{-0.1in}\\
\rdseq{1}{\xarity} &=& \adseq{1}{\xarity} = \bdseq{1}{\xarity} =
\idseq{1}{\xarity}
= \ddseq{1}{\xarity} = 0 &\quad& {\rm\ for\ } \xarity \geq 4; \vspace{-0.1in}\\
\rdseq{n}{0} &=& \adseq{n}{0} = \bdseq{n}{0} = \idseq{n}{0} =
\ddseq{n}{0} = 0 &\quad& {\rm\ for\ }n \geq 1; \vspace{-0.1in}\\
\rdseq{n}{1} &=& \adseq{n}{1} = \bdseq{n}{1} = \idseq{n}{1} =
\ddseq{n}{1} = 0 &\quad& {\rm\ for\ }n \geq 2.
\end{array}
\end{equation*}
%

{}From the production rules
(\mref{prod:rbw}) -- (\mref{prod:unpaired}),
we see that for $n \geq 1, \xarity
\geq 2$:
\begin{equation}
\rdseq{n}{\xarity} = \bdseq{n}{\xarity} + \adseq{n}{\xarity} \mlabel{eq:dbrbw}
\end{equation}
\vspace{-0.35in}
\begin{equation}
\adseq{n}{\xarity} = 2\, \bdseq{n}{\xarity-1} + \bdseq{n}{\xarity-2}
\mlabel{eq:dbassoc}
\end{equation}
\vspace{-0.25in}
\begin{equation}
\bdseq{n}{\xarity} = \idseq{n}{\xarity} + \ddseq{n}{\xarity}
\mlabel{eq:dbround}
\end{equation}
\vspace{-0.25in}
\begin{equation}
\idseq{n}{\xarity} = \ddseq{n-1}{\xarity} + \adseq{n-1}{\xarity}
\mlabel{eq:dbpaired}
\end{equation}

Now for $n \geq 2, \xarity \geq 2$ and any $w \in D(n,\xarity)$, we can
write $w$ uniquely
as $w_{n_1}xw_{n_2} \cdots xw_{n_p}$ where $w_{n_j} \in I(n_j)$ and
$n_1 + \cdots + n_p$ is a composition of $n$ using $p$ positive
integers. Let $\xarity_j$ be the arity of $w_{n_j}$. Then clearly,
$\xarity_1 + \cdots + \xarity_p = \xarity - p + 1$ and so we may refine
Eq.~(\mref{eq:recur}) to:
\begin{equation}
 \ddseq{n}{\xarity} = \sum_{p=2}^{\min(n,\xarity)} \sum_{(\xarity_1, \dots, \xarity_p;\,
\xarity - p+1)} \
\sum_{(n_1,
\dots, n_p;\,n)}
(\idseq{n_1}{ \xarity_1})\cdots (\idseq{n_p}{ \xarity_p}),
\mlabel{eq:dbrecur} \end{equation}
and noting that the case $p=1$ corresponds to a single summand
$i_{n,m}$, refine Eqn.~(\mref{eq:decomp}) to:
\vspace{-0.3in}
\begin{eqnarray}
\hspace{-0.3in}
\bdseq{n}{\xarity} &=& \idseq{n}{\xarity} + \ddseq{n}{\xarity}
=\sum_{p=1}^{\min(n,\xarity)}\! \sum_{(\xarity_1, \dots, \xarity_p;\,
\xarity - p+1)} \
\sum_{(n_1,
\dots, n_p;\,n)}
\!(\idseq{n_1}{ \xarity_1})\cdots (\idseq{n_p}{ \xarity_p})\,.
\mlabel{eq:dbroundrec} \end{eqnarray}
Now from Eq.~(\mref{eq:dbassoc})--(\mref{eq:dbpaired}), we have
\vspace{-0.15in}
\begin{equation*}
\begin{array}{rcl}
\idseq{n}{\xarity} &=& \ddseq{n-1}{\xarity} + \adseq{n-1}{\xarity}\\
&=&\bdseq{n-1}{\xarity} - \idseq{n-1}{\xarity} +
2\,\bdseq{n-1}{\xarity-1}+\bdseq{n-1}{\xarity-2}\,,
\end{array}
\end{equation*}
and so
\begin{equation}
\idseq{n}{\xarity} +\idseq{n-1}{\xarity} = \bdseq{n-1}{\xarity} +
2\,\bdseq{n-1}{\xarity-1} + \bdseq{n-1}{\xarity-2}\,.
\mlabel{eq:dbpairedrec} \end{equation}
Define the bivariate generating series
$${\mathbf R}(z,t) =
\sum_{n=0}^\infty \sum_{\xarity=0}^\infty \rdseq{n}{\xarity}
z^nt^\xarity$$
and similarly define ${\mathbf B}(z,t)$, ${\mathbf I}(z,t), {\mathbf
D}(z,t)$, and ${\mathbf A}(z,t)$.  Note that for  ${\mathbf B}(z,t)$,
${\mathbf I}(z,t), and {\mathbf D}(z,t)$, it does not matter
whether the series indices $n, m$ start at 0 or 1.  Multiply
both sides of Eq.~(\mref{eq:dbpairedrec}) by $z^n t^\xarity$ and
summing up for $n \geq 2, \xarity \geq 2$, the left hand side gives:
\begin{equation*}
\begin{array}{rcl}
&&\displaystyle \sum_{n= 2}^\infty \sum_{\xarity=2}^\infty
\idseq{n}{\xarity} z^n t^\xarity +\sum_{n=
2}^\infty \sum_{\xarity=2}^\infty \idseq{n-1}{\xarity} z^nt^\xarity
\vspace{0.1in}\\
&=& \displaystyle
{\mathbf I}(z,t)-zt+z\sum_{n=
1}^\infty \bigg( - \idseq{n}{1}z^nt +  \sum_{\xarity=1}^\infty
\idseq{n}{\xarity} z^nt^\xarity \bigg)
\vspace{0.1in}\\
&=& \displaystyle
{\mathbf I}(z,t)-zt-z^2t + z\sum_{n=
1}^\infty  \sum_{\xarity=1}^\infty \idseq{n}{\xarity} z^nt^\xarity
\\ \displaystyle
&=&{\mathbf I}(z,t)- zt -z^2t + z\,{\mathbf I}(z,t).
\end{array}
\end{equation*}
Now, we sum the right hand side of Eq.~(\mref{eq:dbpairedrec}) one
term at a time.
\begin{equation*}\begin{array}{rcl}\displaystyle
\sum_{n=2}^\infty \sum_{\xarity=2}^\infty \bdseq{n-1}{\xarity} z^nt^\xarity
&=& \displaystyle
z \sum_{n=1}^\infty \sum_{\xarity=2}^\infty \bdseq{n}{\xarity} z^nt^\xarity
\\
&=& \displaystyle z\,{\mathbf B}(z,t) - z^2t
\end{array}\end{equation*}
\begin{equation*}\begin{array}{rcl}\displaystyle
2\sum_{n=2}^\infty \sum_{\xarity=2}^\infty \bdseq{n-1}{\xarity-1} z^nt^\xarity
&=& \displaystyle
2zt \sum_{n=1}^\infty \sum_{\xarity=1}^\infty \bdseq{n}{\xarity} z^nt^\xarity \\
&=& \displaystyle 2zt\,{\mathbf B}(z,t)
\end{array}\end{equation*}
\begin{equation*}\begin{array}{rcl}\displaystyle
\sum_{n=2}^\infty \sum_{\xarity=2}^\infty \bdseq{n-1}{\xarity-2} z^nt^\xarity &=&
\displaystyle
zt^2 \sum_{n=1}^\infty \sum_{\xarity=0}^\infty \bdseq{n}{\xarity} z^nt^\xarity \\
&=& \displaystyle zt^2 \big({\mathbf B}(z,t)+ \sum_{n=1}^\infty
\bdseq{n}{0}z^n \big)\\ &=& \displaystyle zt^2 \,{\mathbf B}(z,t)
\end{array}\end{equation*}
Hence the right hand side sums to
$z(1+t)^2\, {\mathbf B}(z,t) - z^2t$, giving the identity
\begin{equation}
(1+z)\,{\mathbf I}(z,t) -zt
 = z(1+t)^2 \,{\mathbf B}(z,t)
\mlabel{eq:pairedround}
\end{equation}
Using Eq.~(\mref{eq:dbroundrec}), we have
\allowdisplaybreaks{
\begin{equation*}
\begin{array}{rcl}
&&\displaystyle \sum_{n=1}^\infty \sum_{\xarity=1}^\infty
\bdseq{n}{\xarity} z^nt^\xarity \vspace{0.2in}
\\ &=&\displaystyle  \sum_{n=1}^\infty \sum_{\xarity=1}^\infty
\sum_{p=1}^{\min(n,\xarity)}
\sum_{(\xarity_1, \dots, \xarity_p;\,
\xarity - p+1)} \
\sum_{(n_1,
\dots, n_p;\,n)}
(\idseq{n_1}{ \xarity_1})\cdots (\idseq{n_p}{ \xarity_p})
 z^nt^\xarity \vspace{0.2in} \\
&=&\displaystyle  \sum_{n=1}^\infty \sum_{\xarity=1}^\infty
\sum_{p=1}^{\min(n,\xarity)}
\sum_{(\xarity_1, \dots, \xarity_p;\,
\xarity - p+1)} \
\sum_{(n_1,
\dots, n_p;\,n)}
(\idseq{n_1}{ \xarity_1}z^{n_1}t^{\xarity_1})\cdots (\idseq{n_p}{
\xarity_p}z^{n_p}t^{\xarity_p}) t^{p-1} \vspace{0.2in} \\
&=&\displaystyle \sum_{p=1}^\infty
\bigg( \sum_{\prun=1}^\infty \sum_{\xrun=1}^\infty
\idseq{\prun}{\xrun}z^\prun t^\xrun \bigg)^p t^{p-1},
\end{array}
\end{equation*}
}
and hence
\begin{equation}
{\mathbf B}(z,t)
= \frac{{\mathbf I}(z,t)}{1 - t\, {\mathbf I}(z,t)}.
 \label{eq:dbroundpaired}
\end{equation}
Thus we obtained the identity defining ${\mathbf I}(z,t)$ as
\begin{equation}
(1+z)\,{\mathbf I}(z,t)- zt
= z(1+t)^2\frac{{\mathbf I}(z,t)}{1 - t\, {\mathbf I}(z,t)}
\mlabel{eq:idpaired}
\end{equation}
Solving this quadratic equation in ${\mathbf I}(z,t)$ and using the initial
conditions, we found
\begin{equation}
{\mathbf I}(z,t) = \frac{1-2tz - \sqrt{1-4zt-4zt^2}}{2t(z+1)},
\mlabel{eq:dgfpaired} \vspace{0.1in}\end{equation}
\begin{equation}
{\mathbf B}(z,t) = \frac{1-2zt-2zt^2 -
\sqrt{1-4zt-4zt^2}}{2t(1+t)^2z}, \mlabel{eq:dgfround}
\vspace{0.1in}\end{equation}
\begin{equation}
{\mathbf D}(z,t)\! =\! \frac{2z^2t^3
\!+\!2z^2t^2\!-\!3zt^2\!-\!4zt\!+\!1
\!+\!(zt^!+\!2zt\!-\!1)\sqrt{1\!-\!4zt\!-\!4zt^2}}{2t(1+t)^2z(1+z)}.
\hspace{0.2in} \mlabel{eq:dgfunpaired}
\vspace{0.1in}\end{equation}
We can also obtain the bivariate generating series for
$\adseq{n}{\xarity}$.
\begin{equation*}
\begin{array}{rcl}
{\mathbf A}(z,t) &=& \displaystyle \sum_{n=0}^\infty
\bigg(\adseq{n}{0}z^n + \adseq{n}{1}z^nt+
\sum_{\xarity=2}^\infty \adseq{n}{\xarity}z^nt^\xarity \bigg)
\vspace{0.2in}\\
&=& \displaystyle  t + \sum_{\xarity=2}^\infty \adseq{0}{\xarity} t^\xarity +
\sum_{n=1}^\infty \bigg(
\sum_{\xarity=2}^\infty \adseq{n}{\xarity}z^nt^\xarity \bigg)\vspace{0.2in}\\
&=&\displaystyle  t+
\sum_{n=1}^\infty
\sum_{\xarity=2}^\infty ( 2\,\bdseq{n}{\xarity-1} +
\bdseq{n}{\xarity-2}) z^n t^\xarity \vspace{0.2in}\\
&=&\displaystyle  t+ \sum_{n=1}^\infty \bigg(
t \sum_{\xarity=1}^\infty 2\,\bdseq{n}{\xarity} z^nt^\xarity \bigg) +
t^2 \sum_{n=1}^\infty \sum_{\xarity=0}^\infty \bdseq{n}{\xarity} z^n
t^\xarity \vspace{0.2in} \\ &=& t+ 2t\, {\mathbf B}(z,t) + t^2
\,{\mathbf B}(z,t). \end{array}
\end{equation*}
The bivariate generating function ${\mathbf A}(z,t)$ is given by
\begin{equation}
{\mathbf A}(z,t) = \frac{2+t-2zt-2zt^2
-(2+t)\sqrt{1-4zt-4zt^2}}{2(1+t)^2z}. \mlabel{eq:dgfassoc}
\end{equation}
Interesting sequences and counting information can be derived from
these functions. For example,
the power series expansion for ${\mathbf
A}(z,t)$ in low degrees in $z$, with coefficient in
$t$ accurate up to ${\mathcal O}(t^{16})$, is
\begin{equation*}
t + (2t^2 + t^3)z + (4t^3 + 6t^4+2t^5)z^2 +
(10t^4 + 25t^5 + 20t^6 + 5t^7)z^3 + \cdots\,.
\end{equation*}
{}From this, we can read off that there are 60 associate \RBWs in
$\rbset$ with
three pairs of brackets, and 10 of these have arity 4, 25 have arity
5, 20 have arity 6, and 5 have arity 7.
By expanding the series using $t$ as the main variable, with
coefficient in $z$ accurate up to ${\mathcal{O}}(z^{10})$ we have
\begin{equation*}
t + 2z t + (z+4z^2)t^3 + (6z^2+10z^3)t^4 +
(2z^2 + 25z^3 + 28z^4)t^5 + \cdots\,.
\end{equation*}
{}From this we see that there are 55 associate \RBWs in $\rbset$ with arity 5, and
2 of these have degree 2, 25 have degree 3 and 28 have degree 4. By
specializing $z = 1$, we obtain a sequence for the number of
associates
in $\rbset$ with arity $\xarity$, $(\xarity = 1, 2, ...)$
$$1, 2, 5, 16, 55, 202, 773, 3052, \dots$$
This sequence is new and not in the Sloane data base.

\begin{theorem}
The generating series for $\rdseq{n}{\xarity}$ is
\begin{equation}
{\mathbf R}(z,t)=\frac{1-\sqrt{1-4zt-4zt^2}}{2tz}\,.
\mlabel{eq:dgfrbw}  \end{equation}
\end{theorem}
\proofbegin
This follows by Rule (\mref{prod:rbw}) using
${\mathbf R}(z,t)
=  1 + {\mathbf B}(z,t) + {\mathbf A}(z,t)$ and Eqs.
(\mref{eq:dgfround}) and (\mref{eq:dgfassoc}).
\proofend

Once again we have proved Theorem \mref{thm:zero} (by putting $t =
1$ in Eqs. (\mref{eq:dgfpaired})--(\mref{eq:dgfrbw})).
By specializing $z=1$ to Eq.~(\mref{eq:dgfrbw}), we obtain the
sequence A025227:
$$1,2,4,12,40,144,544,2128, \dots\,,$$
and thus give that sequence a new combinatorial interpretation.
We easily obtain new sequences by specializing to other values, such
as for $z = 2, 3, 4, 5, \dots$:
\begin{equation*}
\begin{array}{rcl}
z = 2:&\quad& 1, 3, 12, 66, 408, 2712, 18912, 136488, \dots \,,\\
z = 3:&& 1, 4, 24, 192, 1728, 16704, 169344, \dots \,,\\
z = 4:&& 1, 5, 40, 420, 4960, 62880, 835840, \dots, \,,\\
z = 5:&& 1, 6, 60, 780, 1140, 178800, 2940000, \dots, \,.
\end{array}
\end{equation*}
Moreover, this result is more refined than Theorem
\mref{thm:zero} and its corollary.
Indeed, we note that Eq.~(\mref{eq:dgfrbw}) is related to the
well-known generating series
\begin{equation}
{\mathbf C}(z) = \sum_{n=0}^\infty C_n z^n = \frac{1 - \sqrt{1-4z}}{2z}
\mlabel{eq:gfcat}
\end{equation}
for the Catalan sequence $C_n$. We have clearly:
\begin{equation*}
\begin{array}{rcl}
{\mathbf R}(z,t) &=& \displaystyle \frac{1-\sqrt{1-4zt(1+t)}}{2tz} =
(1+t)\,{\mathbf C}(zt(1+t)) \\ &=&\displaystyle
\sum_{n=0}^\infty C_n z^nt^n(1+t)^{n+1} \vspace{0.1in}\\
&=&\displaystyle  \sum_{n=0}^\infty  \bigg(\sum_{j=0}^{n+1}
\binomial{n+1}{j} C_n z^n t^{n+j} \bigg) \vspace{0.1in}\\
&=&\displaystyle  \sum_{n=0}^\infty  \bigg(\sum_{\xarity=n}^{2n+1}
\binomial{n+1}{\xarity-n} C_n z^n  t^\xarity \bigg)
\end{array}
\end{equation*}
and hence
\begin{equation}
\rdseq{n}{\xarity} = \begin{cases}\displaystyle {\binomial{n+1}
{\xarity-n}} C_n &{\rm if\ } n \leq \xarity \leq 2n+1, n \geq 0;\cr
\cr
\qquad 0 &{\rm otherwise.}\cr \end{cases}
\end{equation}
This result not only provides the proof anew that $R(n)$ has $2^{n+1}
C_n$
\RBWs, but also that these are distributed by their arities from $n$
to $2n+1$ according to the binomial theorem.
In a similar fashion, using Eq.~(\mref{eq:dgfround})
and a modified
version of Eq.~(\mref{eq:gfcat}):
\begin{equation*}
\overline{{\mathbf C}}(z) = \sum_{n=1}^\infty C_n z^n = \frac{1 -2z -
\sqrt{1-4z}}{2z} \vspace{0.1in}
\end{equation*}
the doubly-indexed sequence for the
number of \round \RBWs with degree $n$ and arity $\xarity$ has the
same property, that is,
\begin{equation}
\bdseq{n}{\xarity} = \begin{cases} \displaystyle {\binomial{n-1}
{\xarity-n}}C_n &{\rm if\ } n \leq \xarity \leq 2n-1, n \geq 1;\cr \cr
\qquad 0 &{\rm otherwise.}\cr  \end{cases}
\end{equation}
This latter distribution, like the one for $R(n,\xarity)$,
 was first observed by experimental computations, but the proof is not
obvious because among the $C_n$ ways to set up the structure of $n$
pairs of balanced brackets, the number of ways to insert $\xarity$ $x$'s to
form \round \RBWs (or \RBWs in the case of $R(n,\xarity)$) depends on the
individual bracket structure (and sometimes,
this number can be zero). An example that illustrates this
observation is the set $B(3,4)$ which has the following 10 elements.
The $C_3 = 5$ possible bracket structures (which correspond to the 5
possible rooted trees with 4 vertices) are shown on the left.

\medskip
\centerline{
\begin{tabular}{|c|c|c|}
\hline
structure & count & \round \RBWs\\
\hline
\hline
$\lc \,\lc \,\lc \,\rc \,\rc \,\rc $&4&$\lc x\lc x\lc x\rc x\rc \rc ,
\lc \lc x\lc x\rc x\rc x\rc , \lc x\lc x\lc x\rc \rc x\rc ,\lc x\lc
\lc x\rc x\rc x\rc $\\ \hline
$\lc \,\lc \,\rc \,\lc \,\rc \,\rc $&2&$\lc x\lc x\rc x\lc x\rc \rc ,
\lc \lc x\rc x\lc x\rc x\rc $\\ \hline
$\lc \,\lc \,\rc \,\rc \,\lc \,\rc $&2&$ \lc x\lc x\rc \rc x\lc x\rc ,
\lc \lc x\rc x\rc x\lc x\rc $\\ \hline
$\lc \,\rc \,\lc \,\lc \,\rc \,\rc $&2&$\lc x\rc x\lc x\lc x\rc \rc ,
\lc x\rc x\lc \lc x\rc x\rc $\\ \hline
$\lc \,\rc \,\lc \,\rc \,\lc \,\rc $&0&\quad\\
\hline
\end{tabular}
}

\subsection{Algorithm for generating $\rbset$}
\mlabel{sec:oneptfour}
It is also interesting to note that the two-variable generating functions
studied above allow us to obtain a new and more effective way to
generate, say $B(n,\xarity)$ recursively, from constructions using \round
\RBWs alone. From Eq.~(\mref{eq:pairedround}) and
Eq.~(\mref{eq:dbroundpaired}), we obtain the identity satisfied by
${\mathbf B}(z,t)$:
\begin{equation}
{\mathbf B}(z,t) - zt =  2zt(1 + t)\,{\mathbf B}(z,t)
+ zt(1+t)^2\, {\mathbf B}(z,t)^2. \mlabel{eq:idround}
\end{equation}
Substituting the definition of ${\mathbf B}(z,t)$ into the above yields:
\begin{equation*}
\begin{array}{rcl}
&&\displaystyle \sum_{n=1}^\infty \sum_{\xarity=1}^\infty
\bdseq{n}{\xarity} z^nt^\xarity - zt \vspace{0.2in}\\
&=&\displaystyle
2zt(1+t)
\bigg(\sum_{n=1}^\infty  \sum_{\xarity=1}^\infty \bdseq{n}{\xarity}
z^nt^\xarity \bigg) +
zt(1+t)^2 \bigg(\sum_{n=1}^\infty  \sum_{\xarity=1}^\infty
\bdseq{n}{\xarity} z^nt^\xarity \bigg)^2
\end{array}
\end{equation*}
from which we obtain, for $(n,\xarity) \neq (1,1)$,
\begin{equation}
\begin{array}{rcl}
\bdseq{n}{\xarity} &=&\displaystyle 2\,\bdseq{n-1}{\xarity-1} +
2\,\bdseq{n-1}{\xarity-2}\ +
\sum_{\tiny\begin{array}{c} (n_1, n_2;\, n-1)\vspace{-0.05in}
\\ (\xarity_1, \xarity_2;\,\xarity-1)\end{array}}
\bdseq{n_1}{\xarity_1}\bdseq{n_2}{ \xarity_2} \\
&\qquad&\displaystyle +\ 2
\sum_{\tiny\begin{array}{c}(n_1, n_2;\,
n-1)\vspace{-0.05in}\\(\xarity_1, \xarity_2;\,\xarity-2)\end{array}}
 \bdseq{n_1}{\xarity_1}\bdseq{n_2}{ \xarity_2}
\ +
\sum_{\tiny\begin{array}{c}(n_1, n_2;\, n-1)
\vspace{-0.05in}\\(\xarity_1, \xarity_2;\,\xarity-3)\end{array}}
\bdseq{n_1}{\xarity_1}\bdseq{n_2}{ \xarity_2}
\end{array} \mlabel{eq:roundnew}
\end{equation}

or, in explicit summation form:
\begin{equation}
\begin{array}{rcl}
\bdseq{n}{\xarity} &=& \displaystyle
2\,\bdseq{n-1}{\xarity-1} + 2\,\bdseq{n-1}{\xarity-2} +
\sum_{\prun=1}^{n-2} \sum_{\xrun=1}^{\xarity-2}
\bdseq{\prun}{\xrun}\bdseq{n-1-\prun}{ \xarity-1-\xrun} \vspace{0.2in}
\\ &\qquad&\displaystyle
+\ 2 \sum_{\prun=1}^{n-2} \sum_{\xrun=1}^{\xarity-3}
 \bdseq{\prun}{\xrun}\bdseq{n-1-\prun}{ \xarity-2-\xrun}
\ +
\sum_{\prun=1}^{n-2} \sum_{\xrun=1}^{\xarity-4}
\bdseq{\prun}{\xrun}\bdseq{n-1-\prun}{ \xarity-3-\xrun}\,.
\end{array} \mlabel{eq:roundexplicit}
\end{equation}

This recursion suggests that we can generate
$B(n,\xarity)$ (and hence also $R(n,\xarity)$) efficiently and irredundantly
from sets $B(\prun, \xrun)$ with $\prun < n, \xrun < \xarity$ in some
manner. This is indeed the case.

\newcommand\step[1]{\noindent {\bf Step} {\rm #1}\\}
\setcounter{figure}{0}

\begin{theorem} \mlabel{thm:two} Given positive integers $n$ and
$\xarity$, the algorithm below returns
the sets $B(n,\xarity)$ \paren{resp.~$I(n,\xarity)$, resp.~$D(n,\xarity)$} of \round
\paren{resp.~\paired, resp.~\unpaired} \RBWs of degree $n$ and arity
$\xarity$. For $n \geq 4$, $B(n,\xarity)$ is expressed
as the disjoint union of sets constructed from $B(\prun,\xrun)$ with
$\prun <
n$ and $\xrun < \xarity$ according to {\rm
Eq.~(\mref{eq:roundexplicit})}. \end{theorem}

\bigskip
\centerline{\bf Algorithm
for \round \RBWs of degree $n$ and arity $\xarity$}

{\rm
\begin{tabbing}
\, \= xx \= xx \= xx \= x \= x \= \kill
\> {\bf Input}: Positive integers $n$, $\xarity$\\
\> {\bf Output}:\\
\> \> (a) the set $B(n,\xarity)$ of \round \RBWs of degree $n$ and
arity $\xarity$,\\
\> \> (b) the set $I(n,\xarity)$ of \paired \RBWs of degree $n$ and
arity $\xarity$,\\
\> \> (c) the set $D(n,\xarity)$ of \unpaired \RBWs of degree $n$ and
arity $\xarity$.\\
\> \step{0. If not($n \leq \xarity \leq 2n-1$), then return three
empty sets.}
\> \> \> \> Generate all \round \paren{resp.~\paired, resp.~\unpaired}
\\
\> \> \> \> \> \RBWs with degree $\prun \leq 3$ and arity $\xrun$
between $\prun$ and $2\prun-1$.\\
\> \> \> \> If $n \leq 3$ then return $B(n,\xarity)$, $I(n,\xarity)$,
$D(n,\xarity)$.\\
\> \step{1. For each \RBW $w \in B(n-1,\xarity-1)$, form two \RBWs}
\> \> \> \> $f_{1,1}(w) = \lc x \,w\,\rc $ and $f_{1,2}(w) = \lc \,w\,x\rc $.\\
\> \step{2. For each \RBW $u \in B(n-1,\xarity-2)$, form the \RBWs
$f_{2}(u) = \lc x\,u\,x\rc $}
\> \step{3. For each $\prun = 1 \dots (n-2)$, each $\xrun=1 \dots
\xarity-2$,
} \> \> \> \>  and each pair of \RBWs $(v, y) \in I(\prun, \xrun)
\times B(n-1-\prun, \xarity-1-\xrun)$,\\
\> \> \> \>  form the \RBW $f_3(v,y) = \lc \,v\,x\,y\,\rc $ \\
\> \step{4. For each $\prun = 1 \dots (n-2)$, each $\xrun=1 \dots
\xarity-2$,
} \> \> \> \>  and each pair of \RBWs $(v, y) \in D(\prun,\xrun)
\times B(n-1-\prun, \xarity-1-\xrun)$\\
\> \> \> \>  form the \RBW $f_{4}(v,y)=\lc \,v\,\rc \,x\,y$\\
\> \step{5. For each \RBW $u \in B(n-1,\xarity-2)$, form the \RBWs
$f_{5}(u)=\lc \,x\,\rc \,x\,u$.}
\> \step{6. For each $\prun = 1 \dots (n-2)$, each $\xrun=1 \dots
\xarity-3$,
} \> \> \> \>  and each pair of \RBWs $(v, y) \in B(\prun,\xrun)
\times B(n-1-\prun, \xarity-2-\xrun)$,\\
\> \> \> \>  form the two \RBWs $f_{6,1}(v,y) = \lc x \,v\rc \,x\,y$ and
$f_{6,2}(v,y) = \lc \,v\,x\rc \,x\,y$.\\
\> \step{7. For each $\prun = 1 \dots (n-2)$, each $\xrun=1 \dots
\xarity-4$,}
\> \> \> \> and each pair of \RBWs $(v, y) \in B(\prun,\xrun) \times
B(n-1-\prun, \xarity-3-\xrun)$,\\
\> \> \> \> form the \RBW $f_7(v,y) = \lc x\,v\,x\rc \,x\,y$.\\
\> \step{8. Return the union of all the \RBWs formed in Steps
1--3 as $I(n,\xarity)$,}
\> \> \> \> the union of all the \RBWs formed in Steps 4--7
as $D(n,\xarity)$,\\
\> \> \> \> and the union of $I(n,\xarity)$ and $D(n,\xarity)$ as $B(n,\xarity)$.
\end{tabbing}
}
\proofbegin By Eq.~(\mref{eq:idround}), we know that $B(n,\xarity)$,
$I(n,\xarity)$, and $D(n,\xarity)$ are empty if $\xarity < n$ or $\xarity  > 2n-1$.
The cases when $n \leq 3$ is taken care of in Step~0.
So suppose $n\geq 4$ and $n \leq \xarity \leq 2n-1$. We first
note the disjoint union $B(n,\xarity)=I(n,\xarity)\cup D(n,\xarity)$
from Eq.~(\mref{prod:round}).

\vspace{-0.15in}
Any word $z$ in $I(n,\xarity)$ is of the form $\lc z'\rc $ where $z'
\in R(n-1,\xarity)$ by Eq.~(\mref{prod:paired}). Then one and exactly one
of the following statements on $z'$ is true.
\begin{list}
{Case--\arabic{enumi}.}{\usecounter{enumi}
\setlength{\rightmargin}{1cm}
\setlength{\leftmargin}{2.2cm}
\setlength{\labelwidth}{1.6cm}
\setlength{\listparindent}{0cm}
\setlength{\labelsep}{0.3cm}
\setlength{\itemsep}{0.2cm}
}
\item \mlabel{case:one} Either $z'$ starts with $x$ but does not end
with $x$, implying $z=\lc x\,w\rc $ with $w\in B(n-1,\xarity-1)$; or\\
$z'$ does not start with $x$ but ends with $x$, implying $z=\lc w\,x\rc $
with $w\in B(n-1,\xarity-1)$;

\item \mlabel{case:two} $z'$ starts with $x$ and ends with $x$,
implying $z=\lc x\,u\,x\rc $ with $u\in B(n-1,\xarity-2)$;

\item \mlabel{case:three} $z'$ neither starts nor ends with $x$. Then
$z'$ must be decomposable and so $z' \in D(n-1,\xarity)$.
Let $v$ be the leftmost \paired subword of $z'$
(see Eq.~(\mref{prod:unpaired})) and let $\prun$ be the degree of $v$
and $\xrun$ be the arity of $v$. Then $1 \leq \prun \leq n-2$, $1
\leq \xrun \leq \xarity-2$, and $v \in I(\prun,\xrun)$. Moreover, we
can
write $z'$ uniquely as $v\,x\,y$ and $z=\lc v\,x\,y\rc $ where $y \in
B(n-1-\prun,\xarity-1-\xrun)$.
\end{list}
This proves that $I(n,\xarity)$ is a subset of the set of \RBWs generated
by Steps 1--3. Conversely, any \RBW $z$ generated by Steps 1--3
clearly belongs to $I(n,\xarity)$ by definition. Thus $I(n,\xarity)$
is precisely the set of \RBWs generated by Steps 1--3, and moreover,
the sets generated in each of these steps are disjoint.
This shows that
\begin{equation} \idseq{n}{\xarity} = 2\,
\bdseq{n-1}{\xarity-1}+\bdseq{n-1}{\xarity-2}+
\sum_{\prun=1}^{n-2} \sum_{\xrun=1}^{\xarity-2}
\idseq{\prun}{\xrun}\,\bdseq{n-1-\prun}{ \xarity-1-\xrun}.
\mlabel{eq:setindecomp}
\end{equation}

Now consider $z \in D(n,\xarity)$. As in Case--\mref{case:three} above,
$z$ is of the form $\lc z'\rc \,x\,y$ for a unique \RBW $z'$ ($\lc z'\rc $ being
the leftmost \paired subword of $z$) and a unique
\round \RBW $y$. Let $\prun$ be the degree of $z'$.
Then $1 \leq \prun \leq n-2$, and one and exactly one of the
following statements on $z'$ is true.

\begin{list}
{Case--\arabic{enumi}.}{\usecounter{enumi}
\setlength{\rightmargin}{1cm}
\setlength{\leftmargin}{2.2cm}
\setlength{\labelwidth}{1.6cm}
\setlength{\listparindent}{0cm}
\setlength{\labelsep}{0.3cm}
\setlength{\itemsep}{0.2cm}
}
\setcounter{enumi}{3}
\item \mlabel{case:four} $z'$ neither starts nor ends with an $x$,
implying that $z'$ is \unpaired. Denoting $z'$ by $v$, and letting
$\xrun$ be the arity of $v$, we have $z=\lc v\rc \,x\,y$ with
$$(v,y)\in
D(\prun,\xrun)\times B(n-1-\prun,\xarity-1-\xrun),$$ and $1\leq
\xrun\leq \xarity-2$.

\item \mlabel{case:five} $z'$ is $x$, implying $z=\lc \,x\,\rc \,x\,y$ and
$y\in B(n-1,\xarity-2)$.

\item  \mlabel{case:six} Either $z'$ starts with $x$ but does not end
with $x$, implying $z=\lc x\,v\,\rc \,x\,y$ for some
$(v,y)\in B(\prun,\xrun)\times B(n-1-\prun,\xarity-2-\xrun)$ where
$1\leq \xrun\leq \xarity-3$; or\\
$z'$ does not start with $x$ but ends with $x$, implying
$z=\lc v\,x\rc \,x\,y$ with $(v,y)\in B(\prun,\xrun)\times
B(n-1-\prun,\xarity-2-\xrun)$ where $1\leq \xrun \leq \xarity-3$;

\item  \mlabel{case:seven} $z'$ starts with $x$ and ends with $x$
(but is not $x$), implying $z=\lc x\,v\,x\rc \,x\,y$ with
$(v,y)\in B(\prun,\xrun)\times B(n-1-\prun,\xarity-3-\xrun)$ where
$1\leq \xrun\leq \xarity-4$;
\end{list}

\noindent
Note that Cases \mref{case:four}--\mref{case:seven}
correspond respectively to Steps 4--7 of the algorithm and generate
disjoint subsets of
$D(n,\xarity)$. This shows that $D(n,\xarity)$ is precisely the set of \RBWs
generated by these steps and hence
\begin{equation}
\begin{array}{rcl}
\ddseq{n}{\xarity} &=& \displaystyle
\sum_{\prun=1}^{n-2} \sum_{\xrun=1}^{\xarity-2}
\ddseq{\prun}{\xrun}\,\bdseq{n-1-\prun}{ \xarity-1-\xrun} +
\bdseq{n-1}{\xarity-2} \vspace{0.2in}\\ &\qquad& \displaystyle
+\ 2 \sum_{\prun=1}^{n-2} \sum_{\xrun=1}^{\xarity-3}
 \bdseq{\prun}{\xrun}\,\bdseq{n-1-\prun}{ \xarity-2-\xrun}
\ +
\sum_{\prun=1}^{n-2} \sum_{\xrun=1}^{\xarity-4}
\bdseq{\prun}{\xrun}\,\bdseq{n-1-\prun}{ \xarity-3-\xrun}\,.
\end{array} \mlabel{eq:setdecomp}
\end{equation}
Combining equations (\mref{eq:setindecomp}) and (\mref{eq:setdecomp})
finishes the proof (and provides a second,
constructive, proof of Eq.~(\mref{eq:roundexplicit})).
\proofend

\section{One generator and one operator: arbitrary exponent case}
\mlabel{sec:two}

We now consider the more general cases.  In this section, we
generalize previous results to the cases of one generator $x$ and one
operator $P=\rbop$ without requiring these to be idempotent.
Referring to Example \mref{ex:four}, we have again $\xnumb~=~1$, but
now $\xexp = \xexp_1$ is arbitrary (including $\xexp_1 = \infty$).
Moreover, we also restrict the set of Rota-Baxter words
$\frakMl(1,(v))$ to those where the number of consecutive applications
of the operator $P$ is bounded by a given $\pexp$ (which may also be
$\infty$, in which case there will be no restriction at all).

\vspace{-0.1in}
\subsection{Notations} \mlabel{sec:twoptone}

\vspace{-0.1in}
For this section, we now introduce new terminology and notations.  For
any \rbw $w$, and operator $P = \rbop$ occurring in $w$, a $P$-{\bf
\run } is any occurrence in $w$ of consecutive compositions of $\rbop$
of maximal length (that is, of immediately nested $\rbop$, where the
length is the number of consecutive applications of $P$).  Recall from
Definition \mref{de:rbgen} for any generator $x$, an $x$-{\bf \run }
is any occurrence in $w$ of consecutive (algebraic) products of $x$ of
maximal length.  We denote a $P$-\run by $P^{(\prunlen)}$ or
$\rbo{\prunlen}$ if $\prunlen$ is its \run length, and an $x$-\run by
$x^\xrunlen$ if $\xrunlen$ is its \run length. When $\prunlen$ or
$\xrunlen$ is 1, we shall often omit the superscript.  Let $\pexp,
\xexp$ be either positive integers or $\infty$ and let
$\rbset_{\pexp,\xexp}$ be the subset of \rbws (including $\bempty$) where the length of
$P$-\runs is $\leq \pexp$ and the length of $x$-\runs is $\leq \xexp$.
These subsets are potential canonical bases of Rota-Baxter algebras on
one generator.  We have seen in Section \mref{sec:zero} that
$\rbset_{\infty, \xexp}$ is the canonical basis of the free
Rota-Baxter algebra $\ncshao(1, (\xexp))$ (Example \mref{ex:four}).
Also $\rbset_{1,1}$ is the canonoical basis of the free Rota-Baxter
algebra with one idempotent generator and one idempotent operator considered in Section \mref{sec:one} (see also \mcite{A-M}).

For convenience, we say the operator $P = \rbop$ has {\bf
\exponent}~$\pexp$ and the generator $x$ has {\bf \exponent}~$\xexp$
if we are enumerating the set $\rbset_{\pexp, \xexp}$.  This would be
the case for Rota-Baxter algebras where the generator $x$ satisfies
$x^{\xexp+1}= x$ and the operator $P$ satisfies $P^{(\pexp + 1)}(w) =
P(w)$ for any $w$.  In this section, our enumeration on \rbws is valid
for any unary operator $P$.  It is not clear under what conditions a
Rota-Baxter operator $P$ would have \exponent $\pexp$ for $\pexp \geq
2$.

For
$\pdeg \geq 1$, let $\rbset_{\pexp,\xexp}(\pdeg)$ be the subset of
$\rbset_{\pexp,\xexp}$ consisting of all \rbws of degree $\pdeg$, and
for $\xarity \geq 1$, let $\rbset_{\pexp,\xexp}(\pdeg, \xarity)$ be
the subset of $\rbset_{\pexp,\xexp}$ consisting of \rbws with degree
$\pdeg$ and arity $\xarity$. Moreover, for $1\leq \prun \leq \pdeg$,
we let $\rbset_{\pexp,\xexp}(\pdeg,\xarity ; \prun)$ be the subset of
$\rbset_{\pexp,\xexp}(\pdeg,\xarity)$ consisting of \rbws where the
$\pdeg$ pairs of balanced brackets are distributed into exactly
$\prun$ $P$-\runs, and for $1\leq \xrun \leq \xarity$, we let
$\rbset_{\pexp,\xexp}(\pdeg,\xarity ; \prun, \xrun)$ be the subset of
$\rbset_{\pexp,\xexp}(\pdeg,\xarity)$ consisting of \rbws where the
$\pdeg$ pairs of balanced brackets are distributed into exactly
$\prun$ $P$-\runs, and the $\xarity$ $x$'s are distributed into
exactly $\xrun$ $x$-\runs.  Except for $\rbset_{\pexp, \xexp}(\pdeg)$,
these subsets are all finite sets, even when
$\pexp, \xexp$ are infinite, and we shall denote their corresponding
cardinalities by replacing $\rbset$ by the lower case $\rbnum$.  Thus,
for example, $\rbnum_{\pexp, \xexp}(\pdeg,\xarity ; \prun, \xrun)$ is
the cardinality for $\rbset_{\pexp, \xexp}(\pdeg,\xarity ; \prun,
\xrun)$ and the count $r_{\pdeg, \xarity}$ of Section
\mref{sec:one} is now denoted by $\rbnum_{1,1}(\pdeg, \xarity)$.  This
convention will
be used for all other (finite) sets of \RBWs we may introduce later.

As an example for the above terms and notations, the \rbw $$w = x^2\lc
x\lc x^3\rc ^{(2)}x^2\rc = xx\lc x\lc \lc xxx\rc \rc xx\rc $$ is an
element in $\rbset_{\pexp,\xexp}(3,8;2,4)$ for any $\pexp\geq 2,
\xexp\geq 3$ since the 3 pairs of balanced brackets occur in 2
$P$-\runs of \run-lengths 1 and 2, and the 8 $x$'s occur in 4
$x$-\runs of \run-lengths 2, 1, 3, 2.

We also define the following generating series:
\begin{equation}
\rbfun_{\pexp,\xexp}(z) = \sum_{\pdeg=0}^\infty
\rbnum_{\pexp,\xexp}(\pdeg)z^\pdeg \mlabel{eq:vuGFn}
\end{equation}
\begin{equation}
\rbfun_{\pexp,\xexp}(z,t) = \sum_{\pdeg=0}^\infty
\sum_{\xarity=0}^\infty \rbnum_{\pexp,\xexp}(\pdeg,\xarity)z^\pdeg
t^\xarity \mlabel{eq:vuGFnm}
\end{equation}
\begin{equation}
\rbfun_{\pexp,\xexp}(z,t; \prunvar) = \sum_{\pdeg=0}^\infty
\sum_{\xarity=0}^\infty
\sum_{\prun=0}^\infty \rbnum_{\pexp,\xexp}(\pdeg,\xarity
;\prun)z^\pdeg t^\xarity \prunvar^\prun \mlabel{eq:vuGFnml}
\end{equation}
\begin{equation}
\rbfun_{\pexp,\xexp}(z,t; \prunvar, \xrunvar) =
\sum_{\pdeg=0}^\infty \sum_{\xarity=0}^\infty
\sum_{\prun=0}^\infty \sum_{\xrun=0}^\infty
\rbnum_{\pexp,\xexp}(\pdeg,\xarity
;\prun, \xrun)z^\pdeg t^\xarity \prunvar^\prun \xrunvar^\xrun
\mlabel{eq:vuGFnmlk}
\end{equation}

We divide our study into three parts. In Section \mref{sec:twopttwo},
we recall some results on compositions. In Section
\mref{sec:twoptthree},
we develop enumeration formulae relating the general cases to the
idempotent case. In the remaining Section \mref{sec:twoptfour},
we compute the
generating functions, and sketch an algorithm for generating the sets
$\rbset_{\pexp,\xexp}(\pdeg,\xarity)$.

\subsection{Compositions of an Integer} \mlabel{sec:twopttwo}

We recall (see, for example,  \mcite{Mac}) a well-known result
on compositions (also called
ordered partitions) of a positive integer $\xarity$. Let
$\comset(\xarity,\xrun,\xexp)$ be the set of
compositions of the integer
$\xarity$ into $\xrun$ positive integer parts, with each part
at most
$\xexp$ and let $\comfun(\xarity,\xrun,\xexp)$ be the size of this
set. When
$\xexp$ is finite, it is easy to see from the multinomial expansion
of $$(\xarityvar + \xarityvar^2 + \cdots + \xarityvar^\xexp)^\xrun$$
by collecting the coefficients of $\xarityvar^\xarity$ that
$$\comfun(\xarity,\xrun,\xexp)=\sum_{
\stackrel{\xrun_1+\cdots+\xrun_\xexp = \xrun}
{\xrun_1+2\xrun_2+\cdots+ \xexp\xrun_\xexp=\xarity}}
{\binomial{\xrun}{\xrun_1, \cdots, \xrun_\xexp}}$$
where $\xrun_j \geq 0$ is the number of times $\xarityvar^j$ is chosen
among the $\xrun$ factors of the product. It also follows that
$$\sum_{\xarity=1}^\infty \comfun(\xarity,\xrun,\xexp)
\xarityvar^\xarity = (\xarityvar + \xarityvar^2 + \cdots +
\xarityvar^\xexp)^\xrun.$$

Thus, for given $\xrun$ and $\xexp$, we find that the left-hand-side
is the generating function $\comGF_{\xrun,\xexp}(\xarityvar)$ for
$\comfun$ with respect to $\xarity$, namely:
\begin{equation}
\comGF_{\xrun,\xexp}(\xarityvar) := \xarityvar^\xrun \left(
\frac{1 - \xarityvar^\xexp}{1-\xarityvar} \right)^\xrun\,.
\mlabel{eq:comGFfin}
\end{equation}

In a similar but simpler way, by expanding
$(\xarityvar + \xarityvar^2 + \cdots)^\xrun$
when $\xexp=\infty$, we see that
\begin{equation}
\comfun(\xarity,\xrun,\infty)={\binomial{\xarity-1}{\xrun-1}}\,,
\mlabel{eq:infpart}
\end{equation}
which is
the number of
compositions of $\xarity$ into $\xrun$ parts, with no
restrictions on the size of each part. We have an associated
generating function
\begin{equation}
\comGF_{\xrun, \infty}(\xarityvar) := \sum_{\xarity = 1}^\infty
\comfun(\xarity,\xrun,\infty) \xarityvar^\xarity =
\sum_{\xarity = 1}^\infty {\binomial{\xarity-1}{\xrun-1}}
\xarityvar^\xarity =
\left( \frac{\xarityvar}{1 - \xarityvar} \right)^\xrun.
\mlabel{eq:comGFinf}
\end{equation}
At the other extreme case, when $\xexp=1$, then
$\comfun(\xarity,\xrun,1)=\delta_{\xarity,\xrun}$ (Kronecker's
$\delta$) with $\comGF_{\xrun,1}(\xarityvar) = \xarityvar^\xrun$.

Regarding the power series ring $\ZZ[[\xarityvar]]$ over the ring
$\ZZ$ of integers as the completion
of the polynomial ring $\ZZ[\xarityvar]$ in the usual sense, we have
$$\comGF_{\xrun, \infty}(\xarityvar) =
\lim_{\xexp\to \infty} \comGF_{\xrun,\xexp}(\xarityvar).$$
This implies that
$$ \comfun(\xarity,\xrun, \infty)=\lim_{\xexp\to \infty}
\comfun(\xarity,\xrun,\xexp).$$
In fact, $\comfun(\xarity,\xrun, \infty)=\comfun(\xarity,\xrun,\xexp)$
for all $\xexp\geq \xarity$.

\begin{lemma} \mlabel{thm:three}
For any power series $\F(\xarityvar)=\sum_{\xrun\geq 1} a_\xrun
\xarityvar^\xrun$ and $1\leq \xexp\leq \infty$, we have
$$ \F \circ \comGF_{1,\xexp}(\xarityvar)
=\sum_{\xarity\geq 1} \big( \sum_{\xrun=1}^\xarity
\comfun(\xarity,\xrun,\xexp) a_\xrun\big) \xarityvar^\xarity.$$
\end{lemma}

\proofbegin
We have
\begin{equation*}
\begin{array}{rcl}
\F \circ \comGF_{1, \xexp}(\xarityvar) &=& \displaystyle
\sum_{\xrun\geq 1} a_\xrun \left (\sum_{j=1}^\xexp
\xarityvar^j \right)^\xrun \\
&=& \displaystyle \sum_{\xrun\geq 1} a_\xrun \comGF_{\xrun, \xexp}(\xarityvar) \\
&=& \displaystyle \sum_{\xrun\geq 1} \sum_{\xarity=1}^\infty a_\xrun
\comfun(\xarity,\xrun,\xexp) \xarityvar^\xarity \\
&=&\displaystyle \sum_{\xarity\geq 1} \big( \sum_{\xrun=1}^\xarity
\comfun(\xarity,\xrun,\xexp) a_\xrun\big) \xarityvar^\xarity,
\end{array}
\end{equation*}
since $\comfun(\xarity,\xrun, \xexp) = 0$ if $\xrun > \xarity$.
$\square$

\subsection{Enumeration} \mlabel{sec:twoptthree}

Recall that in the case when the Rota-Baxter operator $\rbop$ and the
generator $x$ are both idempotent (that is, each of \exponent 1), we
had studied in Section~\mref{sec:one} the set $\rbset_{1,1}(\pdeg) =
R(\pdeg)$, which is the set of all \rbws of degree $\pdeg$ and
the set $\rbset_{1,1}(\pdeg,\xarity)=R(\pdeg,\xarity)$, which is
the set of all \rbws of degree $\pdeg$ and arity $\xarity$.
It was shown in Eq.~(\mref{eq:dgfrbw}) that the doubly indexed sequence
$\rbnum_{\pdeg,\xarity} = \rbnum_{1,1}(\pdeg,\xarity) $ has the
generating function:
$$ \rbfun_{1,1}(z,t):=\rbfun(z,t)=\sum_{\pdeg, \xarity \geq 0}
\rbnum_{1,1}(\pdeg,\xarity) z^\pdeg t^\xarity=
\frac{1-\sqrt{1-4zt-4zt^2}}{2tz}.$$
We have the disjoint union:
\begin{equation}
\rbset_{\pexp,\xexp}(\pdeg, \xarity)=
\coprod_{\prun=1}^\pdeg \coprod_{\xrun=1}^\xarity
\rbset_{\pexp,\xexp}(\pdeg, \xarity; \prun,\xrun).
    \mlabel{eq:rundecomp}
\end{equation}
We note that in order for $\rbset_{1,1}(\pdeg,\xarity;\prun,\xrun)$ to
be non-empty, we must have $\pdeg=\prun$ and $\xarity=\xrun$, in which
case, we have $\rbset_{1,1}(\pdeg,
\xarity;\prun,\xrun)=\rbset_{1,1}(\prun,\xrun)$. For any $w \in
\rbset_{1,1}(\prun, \xrun)$, $w$ has exactly $\prun$ pairs of balanced
brackets, no two of which are immediately nested, and $w$ has exactly
$\xrun$ $x$'s no two of which are next to each other.
We define below what might be called a ``collapsing
map'':
\begin{equation}
\pxscalar:
\rbset_{\pexp, \xexp}(\pdeg, \xarity; \prun, \xrun) \longrightarrow
\rbset_{1,1}(\prun,
\xrun).
\mlabel{eq:collapse}
\end{equation}
Given $w \in \rbset_{\pexp, \xexp}(\pdeg, \xarity; \prun, \xrun)$, we
define $\pxscalar(w)$ to be the \RBW $w_1 \in \rbset_{1,1}(\prun,
\xrun)$ obtained by replacing each of the $\prun$ $P$-\runs that
appears in $w$ by a single
$P$ and each of the $x$-\runs appearing in $w$ by a single
$x$. This map is clearly surjective for
each pair $(\prun, \xrun)$.  We shall refer to $w_1$ as the {\bf
collapse} of $w$.  The map $\pxscalar$ technically depends on $\prun,
\xrun$, but for simplicity, we will omit that in the notation.
Moreover by Eq.~(\mref{eq:rundecomp}), we may by abuse of language,
extend all these maps uniquely to a single one:
$$
\pxscalar:
\rbset_{\pexp, \xexp}(\pdeg, \xarity) \longrightarrow
\coprod_{\prun=1}^\pdeg\, \coprod_{\xrun=1}^\xarity
\rbset_{1,1}(\prun,
\xrun).$$
For each \RBW $w \in \rbset_{\pexp, \xexp}(\pdeg, \xarity; \prun,
\xrun)$, $w$ has exactly $\prun$ $P$-\runs and exactly $\xrun$
$x$-\runs.  For $1 \leq i \leq \prun$, let $\pdeg_i$ be the
\run-length of the $i$-th $P$-\run of $w$ and for $1 \leq j \leq
\xrun$, let $\xarity_j$ be the \run-length of the $j$-th $x$-\run of
$w$.  We have $1 \leq \pdeg_i \leq \pexp$ (since $P$ has
\exponent
$\pexp$) and $1 \leq \xarity_j \leq \xexp$ (since
$x$ has \exponent $\xexp$).  Then to each $w$ there
corresponds the triplet $(w_1, \pdegv, \xarityv)$ where $w_1 =
\pxscalar(w) \in \rbset_{1,1}(\prun, \xrun)$, a composition
$\pdegv = (\pdeg_1, \dots, \pdeg_\prun)$ of $\pdeg$ into $\prun$ parts
with each part $\pdeg_i$ no bigger than $\pexp$, and a composition
$\xarityv = (\xarity_1, \dots, \xarity_\xrun)$ of $\xarity$
into $\xrun$ parts with no part $\xarity_j$ bigger than $\xexp$.
Conversely, every such triplet determines a unique $w \in
\rbset_{\pexp, \xexp}(\pdeg, \xarity; \prun, \xrun)$ such that
$\pxscalar(w)=w_1$.  Using the sets $\comset(\xarity, \xrun,
\xexp)$ for compositions defined in Section \mref{sec:twopttwo}, we
have a bijection: \begin{equation}
 \rbset_{\pexp, \xexp}(\pdeg, \xarity; \prun, \xrun)
\longleftrightarrow
\rbset_{1,1}(\prun, \xrun) \times  \comset(\pdeg, \prun, \pexp) \times
\comset(\xarity, \xrun, \xexp) \mlabel{eq:collapsebij}
\end{equation}

\begin{theorem} \mlabel{thm:four}
With the  notations established in {\rm Sections~\mref{sec:twoptone}},
{\rm \mref{sec:twopttwo}} and above, the numbers $\rbnum_{\pexp,
\xexp}(\pdeg, \xarity)$ of
\RBWs of degree $\pdeg$ and arity $\xarity$ in the set
$\rbset_{\pexp, \xexp}$ of \rbws with one generator $x$ of exponent
$\xexp$ and one operator $P$ of exponent $\pexp$ are given by:

\renewcommand{\theenumi}{\arabic{enumi}}  
\renewcommand{\theenumii}{\alph{enumii}}  
\renewcommand{\labelenumi}{{\rm (\arabic{enumi})}}  
\renewcommand{\labelenumii}{{\rm (\alph{enumii})}}  

\begin{enumerate}
\item For $1 \leq \pexp \leq \infty$,
$\displaystyle \rbnum_{\pexp, 1}(\pdeg, \xarity)=\sum_{\prun=1}^\pdeg
\comfun(\pdeg,\prun,\pexp) \rbnum_{1,1}(\prun,\xarity)$;
\vspace{0.05in}

\item For $1 \leq \xexp \leq \infty$,
$\displaystyle \rbnum_{1, \xexp}(\pdeg,
\xarity)=\sum_{\xrun=1}^\xarity
\comfun(\xarity,\xrun,\xexp) \rbnum_{1,1}(\pdeg,\xrun)$;
\vspace{0.1in}

\item For $1 \leq \pexp, \xexp \leq
\infty$, $\rbnum_{\pexp,\xexp}(\pdeg,
\xarity)$ is given by any of the three expressions:
\vspace{0.05in}
\begin{enumerate}
\item $\displaystyle \sum_{\xrun=1}^\xarity
\comfun(\xarity,\xrun,\xexp) \rbnum_{\pexp,1}(\pdeg,\xrun)$,
\vspace{0.05in}
\item $\displaystyle \sum_{\prun=1}^\pdeg \comfun(\pdeg,\prun,\pexp)
\rbnum_{1,\xexp}(\prun, \xarity)$,
\vspace{0.05in}
\item $\displaystyle \sum_{\prun=1}^\pdeg \sum_{\xrun=1}^\xarity
\comfun(\pdeg,\prun,\pexp) \comfun(\xarity,\xrun,\xexp)
 \rbnum_{1,1}(\prun,\xrun)$.
\end{enumerate} \end{enumerate}
\end{theorem}
\proofbegin
The proofs are all similar. For example, 3(c) follows from the
bijection in Eq.~(\mref{eq:collapsebij}) and the disjoint union in
Eq.~(\mref{eq:rundecomp}). \proofend

\begin{coro} \mlabel{thm:five}
Suppose one or both of $u$ and $v$ are $\infty$. Then we have
\begin{equation*}
\begin{array}{rcl}
\rbnum_{\infty, 1}(\pdeg, \xarity) &=& \displaystyle
\sum_{\prun=1}^\pdeg {\binomial{\pdeg-1}{\prun-1}}
{\binomial{\prun+1}{\xarity -
\prun}}C_\prun, \vspace{0.2in}\\
\rbnum_{1, \infty}(\pdeg, \xarity) &=& \displaystyle
\sum_{\xrun=1}^\xarity
{\binomial{\xarity-1}{\xrun-1}} {\binomial{\xrun+1}{\pdeg -
\xrun}}C_\prun \vspace{0.2in}\\
\rbnum_{\infty, \infty}(\pdeg, \xarity) &=& \displaystyle   \sum_{\prun=1}^\pdeg
\sum_{\xrun=1}^\xarity {\binomial{\pdeg-1}{\prun-1}}
{\binomial{\xarity-1}{\xrun-1}} {\binomial{\prun+1}{\xrun -
\prun}}C_\prun.
\end{array}
\end{equation*}
\end{coro}

\subsection{Generating Functions and Enumeration Algorithm}
\mlabel{sec:twoptfour}

We now find the generating functions $\rbfun_{\pexp, \xexp}(\pdegvar,
\xarityvar)$ of the number sequences
$\rbnum_{\pexp,\xexp}(\pdeg, \xarity)$ for $1\leq \pexp,\xexp\leq
\infty$ (see Eqs.~(\mref{eq:vuGFnm})--(\mref{eq:vuGFnmlk}) for
definitions and notations).
Recall from Eqs.~(\mref{eq:dgfrbw})--(\mref{eq:gfcat}), we have
$$\rbfun_{1,1}(\pdegvar,\xarityvar)=\rbfun(z,t) =
\frac{1-\sqrt{1-4 \pdegvar \xarityvar (1+\xarityvar)}}{2 \pdegvar
\xarityvar} = (1+\xarityvar) {\mathbf C}(\pdegvar \xarityvar(1+ \xarityvar))$$
where
$$\displaystyle {\mathbf C}(\pdegvar)=\sum_{\pdeg=0}^\infty C_\pdeg
\pdegvar^\pdeg =\frac{1-\sqrt{1-4 \pdegvar}}{2 \pdegvar}$$
is the generating function of the Catalan numbers.

\begin{theorem} \mlabel{thm:six}
Let $1\leq \pexp,\xexp\leq \infty$. The generating function
$\rbfun_{\pexp, \xexp}(\pdegvar, \xarityvar)$ for the number
$\rbnum_{\pexp, \xexp}(\pdeg, \xarity)$ of \RBWs in a Rota-Baxter
algebra with one
operator $P$ with exponent $\pexp$ and one generator $x$ with exponent
$\xexp$ is given by:

\renewcommand{\theenumi}{\arabic{enumi}}  
\renewcommand{\theenumii}{\alph{enumii}}  
\renewcommand{\labelenumi}{{\rm (\arabic{enumi})}}  
\renewcommand{\labelenumii}{{\rm (\alph{enumii})}}  

\begin{enumerate}
\item $\displaystyle \rbfun_{\pexp,
1}(\pdegvar,\xarityvar)=\rbfun_{1,1}\left
(\comGF_{1,\pexp}( \pdegvar), \ \xarityvar \right).$

\item
$\displaystyle \rbfun_{1, \xexp}(\pdegvar, \xarityvar)=\rbfun_{1,1}
\left (\pdegvar, \comGF_{1,\xexp}( \xarityvar) \right).$

\item
$\displaystyle \rbfun_{\pexp,\xexp}(\pdegvar,\xarityvar)=\rbfun_{1,1}
\left (\comGF_{1,\pexp}( \pdegvar), \comGF_{1,\xexp}( \xarityvar)\right).$
\end{enumerate}
where $\comGF_{1,\pexp}$ and $\comGF_{1,\xexp}$ are given by {\rm Eq.}
\paren{{\rm \mref{eq:comGFfin}}} \paren{for finite $\pexp, \xexp$} and by
{\rm Eq.} \paren{{\rm \mref{eq:comGFinf}}}
\paren{for infinite $\pexp, \xexp$}.
\end{theorem}
\proofbegin
The proofs follow from Theorem~\mref{thm:four} and
Lemma~\mref{thm:three}. We just prove case (2):
\begin{equation*}
\begin{array}{rcl}
\rbfun_{1, \xexp}(\pdegvar,\xarityvar)&=&
\displaystyle \sum_{\pdeg,\xarity
\geq 1} \rbnum_{1, \xexp}(\pdeg, \xarity) \pdegvar^\pdeg
\xarityvar^\xarity\\
&=& \displaystyle \sum_{\pdeg, \xarity\geq 1} \sum_{\xrun=1}^\xarity
\comfun(\xarity,\xrun,\xexp) \rbnum_{1,1}(\pdeg,\xrun)
   \pdegvar^\pdeg \xarityvar^\xarity\\
&=& \displaystyle \sum_{\xarity\geq 1}\left( \sum_{\xrun=1}^\xarity
\comfun(\xarity,\xrun,\xexp) \left(\sum_{\pdeg \geq 1}
\rbnum_{1,1}(\pdeg, \xrun)\pdegvar^\pdeg \right)
\right) \xarityvar^\xarity\\
&=& \displaystyle \sum_{\xrun\geq 1} \left(\sum_{\pdeg \geq 1}
\rbnum_{1,1}(\pdeg,\xrun) \pdegvar^\pdeg     \right )
\left(\comGF_{1,\xexp}(\xarityvar) \right)^\xrun\\
&=& \rbfun_{1,1} \left(\pdegvar,\
\comGF_{1,\xexp}(\xarityvar)\right).\  \square
\end{array}
\end{equation*}

\begin{coro} \mlabel{thm:seven}
Suppose one or both of $\pexp, \xexp$ are $\infty$, then the
generating
function for $\rbnum_{\pexp, \xexp}(\pdeg, \xarity)$ is given
by
\begin{equation*}
\begin{array}{rcl}
\rbfun_{\infty, 1}(\pdegvar, \xarityvar) &=&  \displaystyle
\rbfun_{1,1}\left(\frac{\pdegvar}{1-\pdegvar}, \xarityvar \right),\\
\rbfun_{1, \infty}(\pdegvar, \xarityvar) &=&  \displaystyle
\rbfun_{1,1}\left(\pdegvar, \frac{\xarityvar}
{1-\xarityvar}\right),\\
\rbfun_{\infty, \infty}(\pdegvar, \xarityvar) &=& \displaystyle
\rbfun_{1,1}\left(\frac{\pdegvar}{1-\pdegvar}, \frac{\xarityvar}
{1-\xarityvar}\right).
\end{array}
\end{equation*}
\end{coro}

We end this section with a brief description for an algorithm to
enumerate the sets $\rbset_{\pexp, \xexp}(\pdeg, \xarity)$.  The
details will be left out since by means of the disjoint union in
Eq.~(\mref{eq:rundecomp}) and the bijection in
Eq.~(\mref{eq:collapsebij}), this is fairly straight forward.  We
already have an algorithm (see Theorem \mref{thm:two}) for the
enumeration of $\rbset_{1,1}(\prun, \xrun)$ for any positive $\prun,
\xrun$.  We need an algorithm to generate all the compositions of
$\xarity$ in $\comset(\xarity, \xrun, \xexp)$, which would of course
generate $\comset(\pdeg, \prun, \pexp)$, too.  Now the set of
compositions $\xarityv = (\xarity_1, \dots, \xarity_\xrun)$ of
$\xarity$ into exactly $\xrun$ parts without restrictions on the parts
can be enumerated by readily available, efficient, and well-known
algorithms (see COMP\_NEXT of SUBSET library in \mcite{wilf} for
example).  Those compositions whose parts violate the restrictions
$\xarity_j \leq \xexp$ can be easily discarded by modifying the code.

\section{The general case: multiple generators and operators}
\mlabel{sec:three}

In this section, we consider sets of \rbws with $\pnumb$ unary
operators $P_1, \dots, P_\pnumb$ and $\xnumb$ generators $x_1, \dots,
x_\xnumb$.
Since our purpose here is enumeration, we will for now consider formal
bracketed words with brackets ${\lc_i}\,{\rc_i}$
corresponding to
$P_i$ $(1\leq i \leq \pnumb)$ and ignore any other properties of the
operators.

While it is possible to give a definition generalizing that of
Definition \mref{de:apwao}, it is more straight-forward to generalize
the grammar of
Section \mref{sec:one}. Since we do not need such details, we leave
the specification of the grammar to the reader.

We adopt the convention that a vector quantity using the
same symbol as the corresponding scalar quantity will have components
with the same symbol but subscripted. For example, $\pv = (P_1,
\dots, P_\pnumb)$ and $\xv = (x_1, \dots, x_\xnumb)$. The
$\pv$-\exponent vector $\pexpv = (\pexp_1, \dots, \pexp_\pnumb)$
will mean that the operator $P_i$ has \exponent $\pexp_i$, $1 \leq i
\leq \pnumb$, and this means that we only consider \rbws in which the
number of consecutive
applications of the operator $P_i$ is bounded by $\pexp_i$ for each
$i$.  We define $P_i$-\runs (resp.~$x_j$-\runs) similarly to $P$-\runs
(resp.~$x$-\runs),
treating each $P_i$ (resp.~$x_j$) as single operator
(resp.~generator). We shall call a \run of $P$'s (with whatever subscripts) a
$\pv$-\run, and similarly define an $\xv$-\run. Sets of \RBWs
are defined also with the parameters vectorized. For a numerical
vector such as $\pdegv$, we let $|\pdegv|$ denote its {\bf norm}, which is the sum
of its components. Thus, without further explanation, $
|\xarityv|$ will be the total $\xv$-arity of a \RBW $w$ whose
$\xv$-arity vector is $\xarityv$. As
an example, when $\pnumb = \xnumb = 2$, the \RBW
\begin{equation}
\begin{array}{rcl}
w &=& x_1^3 x_2^4 P_1
P_2^{(3)}(x_1 x_2 P_1(x_1))\\
&=& x_1 x_1 x_1 x_2 x_2 x_2 x_2 \lc _1 \lc _2 \lc _2 \lc _2 x_1 x_2
\lc _1 x_1 \rc _1 \rc _2 \rc _2 \rc _2 \rc _1
\end{array}
\mlabel{eq:exrbw-m}
\end{equation}
has three $x_1$-\runs of lengths 3, 1, and 1; two $x_2$-\runs of
lengths 4 and 1; three $\xv$-runs of lengths 7, 2, and 1; two
$P_1$-\runs of lengths 1 and 1; one $P_2$-\run of length 3, and two
$\pv$-\runs of lengths 4 and 1. The $\pv$-degree vector of $w$ is
$(2, 3)$, its $\pv$-degree is 5, its $\xv$-degree vector is $(5, 5)$,
and its $\xv$-arity is 10.

For any given positive integers $\pnumb, \xnumb$, and corresponding
$\pv$-\exponent vector $\pexpv$, $\xv$-\exponent vector $\xexpv$, let
$\rbset_{\pexpv, \xexpv}$ denote the set of \rbws with $\pnumb$
operators and $\xnumb$ generators with corresponding \exponents
vectors $\pexpv$ and $\xexpv$ respectively. The values $\pnumb,
\xnumb$
are implicitly given by the dimensions of the vectors $\pexpv$ and
$\xexpv$
respectively.  In particular $\rbset_{\infty, \infty}$ is the set of
all \rbws in the setting with one generator (of \exponent $\infty$)
and one
operator (also of \exponent $\infty$).  Thus $\rbset_{\infty, \infty}$
is the canonical
basis for $\ncshao(\bfk[x])$ (see Example \mref{ex:unipoly}).
 Without further explanations,
the notations established in Section \mref{sec:twoptone}
will be generalized to their vectorized versions in the obvious way.
For example, $\rbset_{\pexpv, \xexpv}(\pdeg, \xarity;\prun, \xrun)$
will denote the set of all \RBWs of $\rbset_{\pexpv, \xexpv}$ with
$\pv$-degree $\pdeg$ distributed into exactly $\prun$ $\pv$-\runs, and
$\xv$-arity $\xarity$ distributed into exactly $\xrun$ $\xv$-\runs.
Furthermore, the cardinality of this set is denoted by
$\rbnum_{\pexpv, \xexpv}(\pdeg, \xarity; \prun, \xrun)$ and analogous
to Eq.~(\mref{eq:rundecomp}), we have the disjoint union:
\begin{equation}
\rbset_{\pexpv,\xexpv}(\pdeg, \xarity)=
\coprod_{\prun=1}^\pdeg \coprod_{\xrun=1}^\xarity
\rbset_{\pexpv,\xexpv}(\pdeg, \xarity; \prun,\xrun).
    \mlabel{eq:vrundecomp}
\end{equation}
\subsection{The forgetful maps and coloring} \mlabel{sec:threeptone}

We define another family of maps to relate the
general case to the single operator and single generator case we have
solved.  Basically, the forgetful maps removes all the subscripts on
the operators and generators.  By composing the forgetful maps with
the collapsing maps, we further relate the general case with the
case when the single operator and single generator have exponent 1.
These maps capture the structure of \runs of \RBWs in the general
case.  We now set up the notations to make these precise.

 We define the
``forgetful'' maps:
\begin{equation}
\pxmapscalar:  \rbset_{\pexpv, \xexpv}(\pdeg, \xarity; \prun, \xrun)
\longrightarrow \rbset_{\infty, \infty}(\pdeg, \xarity; \prun, \xrun)
\mlabel{eq:pxmapscalar}
\end{equation}
by defining $\pxmapscalar(w)$, for each $w \in \rbset_{\pexpv,
\xexpv}(\pdeg, \xarity; \prun, \xrun)$, to be the \RBW $w_\infty$ in
$\rbset_{\infty, \infty}$ when every $P_i = \lc _i\, \rc _i$, $1\leq
i\leq \pnumb$ that
appears in $w$ is replaced by $P=\rbop$ and every $x_j, 1\leq j\leq
\xnumb$ that appears in $w$ is replaced by $x$. For example, for the
word $w$ in Eq. (\mref{eq:exrbw-m}), we have
$$
\pxmapscalar(w) = x^7 P^{(4)}(x^2 P(x))
= x^7 \lc  x^2 \lc x \rc \rc^{(4)}.
$$
The map $\pxmapscalar$ is clearly surjective and depends on $\pnumb,
\xnumb, \pexpv, \xexpv, \pdeg, \xarity, \prun, \xrun$, but for
notational brevity, we will not explicitly mention $\pnumb, \xnumb,
\prun$, and $\xrun$. Again, by abuse of language and the disjoint
unions in Eq.~(\mref{eq:vrundecomp}) and Eq.~(\mref{eq:rundecomp}), we
may extend $\pxmapscalar$ uniquely to a surjection:
$$\pxmapscalar: \rbset_{\pexpv, \xexpv}(\pdeg, \xarity)
\longrightarrow \rbset_{\infty, \infty}(\pdeg, \xarity).$$
We may further extend $\pxmapscalar$ to a surjection:
$$\pxmap_{\pexpv,\xexpv} : \rbset_{\pexpv, \xexpv}
\longrightarrow \rbset_{\infty, \infty}.$$

It is convenient to refer to the \RBW $w_\infty = \pxmap_{\pexpv,
\xexpv}(w)$ as the {\bf monochrome image} of $w$ and say $w$ is a {\bf
coloring} of $w_\infty$.  A special case of particular importance is
the following.

Consider an $\xv$-run $w$ which we may write in the form
$x_{j_1}^{\xb_1} \cdots x_{j_\xindex}^{\xb_\xindex}$ where $\xindex$
is some positive integer and for $1 \leq d \leq \xindex$, we have $1
\leq j_d \leq \xnumb$, $1 \leq \xb_d \leq \xexp_{j_d}$, and if $d <
\xindex$, then $j_d \ne j_{d+1}$.  The monochrome image of $w$ will be
$x^\xb$, where $\xb = \xb_1 + \cdots + \xb_\xindex$.  Conversely, for
any positive integer $\xb$, any coloring of $x^\xb$ is an
$\xv$-run $w$ of this form.  Thus there is a bijection $\sigma =
\sigma_{\xb, \xnumb, \xexpv}$ between the set of $\xv$-runs
$w$ whose monochrome image is $x^\xb$ and the set
$\colorset(\xb, \xnumb, \xexpv)$ of colorings of $\xb$ identical
objects in a row using up to $\xnumb$ colors, say colors $1, 2, \dots,
\xnumb$, repetitions allowed, so that for $1 \leq j \leq \xnumb$, any
\run of color $j$ has length no longer than $\xexp_j$.  Explicitly,
given $w$, we define a coloring $\sigma(w) \in \colorset(\xb, \xnumb,
\xexpv)$ of the $\xb$ objects using color $j_1$ for the first $\xb_1$
objects, then color $j_2$ for the next $\xb_2$ objects, etc.  The
number of such colorings is denoted by $\colornum(\xb, \xnumb,
\xexpv)$.

Similarly, for any positive integer $\pa$, we have a bijection $\pi =
\pi_{\pa, \pnumb, \pexpv}$ between the set of $\pv$-runs whose
monochrome image is $P^{(\pa)}$ and the set of colorings
$\colorset(\pa,\pnumb, \pexpv)$.  We can also apply these coloring
maps to $\xv$-runs and $\pv$-runs occurring in
any \RBW $w \in \rbset_{\pexpv, \xexpv}(\pdeg,
\xarity;\prun, \xrun)$ whose monochrome image is a fixed $w_\infty =
\pxmapscalar(w)$. More explicitly,  given $w_\infty \in
\rbset_{\infty, \infty}(\pdeg, \xarity; \prun, \xrun)$,
let $(w_1, \pdegv, \xarityv)$ be the triplet corresponding to
$w_\infty$ according to the bijection from Eq.~(\mref{eq:collapsebij})
for the case $\pexp = \xexp = \infty$, so that
$\pdegv = (\pdeg_1, \dots, \pdeg_\prun)$ is a composition of $\pdeg$
into $\prun$ parts,  $\xarityv = (\xarity_1, \dots, \xarity_\xrun)$
is a composition of $\xarity$ into $\xrun$ parts, and, in particular,
$w_\infty$ collapses to $w_1 \in \rbset_{1,1}(\prun, \xrun)$ (that
is, $\pxinfty(w_\infty) = w_1$).    This implies
that for all $w \in
(\pxmap_{\pexpv,\xexpv})^{-1}(w_\infty)
= (\pxmapscalar)^{-1}(w_\infty)$, $w$ has
exactly $\prun$ $\pv$-\runs of lengths $\pdeg_1, \dots, \pdeg_\prun$,
providing the colorings via the maps $\pi_{\pdeg_i, \pnumb, \pexpv}$
for the corresponding $P$-\runs in the monochromatic image $w_\infty$
and similarly $w$ has exactly $\xrun$ $\xv$-\runs of lengths
$\xarity_1, \dots, \xarity_\xrun$, providing the colorings via the
maps $\sigma_{\xarity_j, \xnumb, \xexpv}$ for the corresponding
$x$-runs in $w_\infty$.
This yields the bijection:
\begin{equation}
\rho_{w_\infty}: (\pxmapscalar)^{-1}(w_\infty) \longrightarrow
\left(\prod_{i=1}^\prun \colorset(\pdeg_i, \pnumb,
\pexpv)\right) \times
\left( \prod_{j=1}^\xrun \colorset(\xarity_j, \xnumb,
\xexpv)\right)\,.
\mlabel{eq:invcolor}
\end{equation}
By Eq.~(\mref{eq:collapsebij}), as $w_\infty$ runs through the set
$\rbset_{\infty, \infty}(\pdeg,
\xarity; \prun, \xrun)$, its triplet $(w_1, \pdegv, \xarityv)$ runs
through $\rbset_{1,1}(\prun, \xrun) \times \comset(\pdeg, \prun,
\infty) \times \comset(\xarity, \xrun, \infty)$. Thus we have a
bijection:
\begin{equation}
\begin{array}{rcl}
&&\hspace{-0.2in} \displaystyle \rho:
\rbset_{\pexpv,
\xexpv}(\pdeg,\xarity;\prun,\xrun)
= (\pxmapscalar)^{-1}(\rbset_{\infty,
\infty}(\pdeg, \xarity;\prun,\xrun)) \longleftrightarrow \\
&&\hspace{-0.2in} \displaystyle  \ \rbset_{1,1}(\prun, \xrun)
\times \coprod_{\pdegv\,\in \comset(\pdeg, \prun, \infty)}
\left(\prod_{i=1}^\prun \colorset(\pdeg_i, \pnumb,
\pexpv)\right) \times \coprod_{\xarityv\,\in
\comset(\xarity,\xrun,\infty)}
\left( \prod_{j=1}^\xrun \colorset(\xarity_j, \xnumb,
\xexpv)\right)
\end{array}\mlabel{eq:colcol}
\end{equation}
We gather below two results on coloring.

\begin{lemma} Let $\xexpv = (\xrunlen, \dots, \xrunlen)
$ be a vector with $\xnumb$ identical coordinates where $\xrunlen$ is
either an integer $\geq 1$ or $\infty$. Then  \begin{equation}
\colornum(\xb, \xnumb, \xexpv)
 = \begin{cases}
1 &{\rm \ if\ } q = 1 {\rm\ and\ } b \leq \xrunlen; \cr
0 &{\rm \ if\ } q = 1 {\rm\ and\ } b > \xrunlen; \cr
\displaystyle \sum_{\xindex=1}^\xb \comfun(\xb,
\xindex, \xrunlen) \xnumb (\xnumb -1)^{\xindex - 1} &{\rm\ if\ } q
\geq 2. \cr
\end{cases}
\mlabel{eq:colorcount}
\end{equation}
If $\xrunlen$ is infinite, then this simplifies to:
\begin{equation}
\colornum(\xb, \xnumb, \inftyv) =  \xnumb^\xb.
\mlabel{eq:colorinf}
\end{equation}
\mlabel{thm:lemma8}
\end{lemma}

\vspace{-0.2in}
\proofbegin The case when $q=1$ is obvious.
Suppose $q \geq 2$ and $\xrunlen$ is finite (resp.~infinite).  The
number $\xindex$
of same color \runs in a coloring of $\xb$ objects can be at most
$\xb$.  For each composition $\xb_1 + \cdot + \xb_\xindex = \xb$
in $\comset(\xb, \xindex, \xrunlen)$, the \run-lengths $\xb_1, \dots,
\xb_\xindex$, which are uniformly bounded by $\xrunlen$
(resp.~unbounded), are fixed (therefore, the locations of the
$\xindex$ \runs are also fixed).  The uniform bound
(resp.~unboundedness) allows any \run be assigned any color, except
that adjacent \runs must have different colors.  Thus there are
$\xnumb (\xnumb-1)^{\xindex -1}$ ways to choose the colors.  Any two
such colorings with distinct parameters are distinct.  This proves the
first statement.

Now suppose $\xrunlen$ is infinite. Then we can color any of the $\xb$
objects with any of the $\xnumb$ colors and there are $\xnumb^\xb$
ways. Alternatively, the second statement also follows from the first
using Eq.~(\mref{eq:infpart}) and the Binomial Theorem.
\proofend

\begin{coro} \mlabel{thm:corol9} Let $\xexpv = (\xrunlen, \dots,
\xrunlen)$ be the same as in {\rm Lemma \mref{thm:lemma8}}.  Then
the generating series for $\colornum(\xb, \xnumb, \xexpv)$ for fixed
$\xnumb$ and $\xexpv$ is given by
$$
\colorGF_{\xnumb, \xexpv}(\xarityvar) =
\begin{cases} \displaystyle
\comGF_{1, \xrunlen}(\xarityvar) &{\rm\ if\ } q = 1; \vspace{0.2in}\cr
\displaystyle
\frac{\xnumb\, \comGF_{1, \xrunlen}(\xarityvar)}
{1 - (\xnumb -1) \comGF_{1,\xrunlen}(\xarityvar)} &{\rm\ if\ }q \geq
2. \end{cases}$$
If $\xrunlen = \infty$, then this simplifies to
$$\colorGF_{\xnumb, \inftyv}(\xarityvar) = \frac{\xnumb \xarityvar}{1
- \xnumb \xarityvar}.$$ \end{coro}

\proofbegin
The case $q=1$ is clear. For $q \geq 2$, the generating
function is
\begin{equation*}
\colorGF_{\xnumb, \xexpv}(\xarityvar) =
\sum_{\xb=1}^\infty \left(\sum_{\xindex=1}^\xb \comfun(\xb,
\xindex, \xrunlen) \xnumb (\xnumb -1)^{\xindex - 1}\right)
\xarityvar^\xb
\end{equation*}
which, by Lemma \mref{thm:three} applied to
$$\F(\xarityvar)
= \sum_{\xindex=1}^\infty \xnumb (\xnumb-1)^{\xindex-1}
\xarityvar^\xindex = \frac{\xnumb \xarityvar}{ 1 -
(\xnumb-1)\xarityvar}\,,$$
is $\F(\comGF_{1, \xrunlen}(\xarityvar))$, as required. The $\xrunlen
= \infty$ case follows from Eqs.~(\mref{eq:colorinf}) and
(\mref{eq:comGFinf}).  \proofend

\vspace{-0.4in}
We say the operators (resp.~generators) have {\bf uniform \exponents}
$\prunlen$ (resp.~$\xrunlen$) if
all the coordinates of the $\pv$-\exponent
(resp.~$\xv$-\exponent) vector $\pexp$ (resp.~$\xexp$)
are equal to $\prunlen$ (resp.~$\xrunlen$), in which case we write
$\pexpv = \vec{\prunlen}$ (resp.~$\xexpv = \vec{\xrunlen}$).  In the
next two subsections, we  consider the infinite
\exponent case, followed by the finite uniform \exponent case.

\subsection{
Uniform \exponents: Infinite cases}
\mlabel{sec:threepttwo}

\vspace{-0.1in}
\begin{theorem} \mlabel{thm:eight}
Let $\rbset_{\pexpv, \xexpv}$ be the set of Rota-Baxter words
with $\pnumb$ operators $\pv$ and $\xnumb$ generators $\xv$, with
corresponding uniform \exponent vectors $\pexpv$ and
$\xexpv$ respectively. Then for
any positive integers $\pdeg, \xarity$, we have
\vspace{-0.1in}
\renewcommand{\theenumi}{\arabic{enumi}}  
\renewcommand{\theenumii}{\alph{enumii}}  
\renewcommand{\labelenumi}{{\rm (\arabic{enumi})}}  
\renewcommand{\labelenumii}{{\rm (\alph{enumii})}}  
\begin{enumerate}
\item $\rbnum_{\onev,\inftyv}(\pdeg,\xarity)
    =\xnumb^\xarity\, \rbnum_{1,\infty}
(\pdeg,\xarity)$,
\item $\rbnum_{\inftyv,\onev}(\pdeg,\xarity)
    =\pnumb^\pdeg\, \rbnum_{\infty,1}
(\pdeg,\xarity)$,
\item $\rbnum_{\inftyv,\inftyv}(\pdeg,\xarity)
    =\pnumb^\pdeg\, \xnumb^\xarity\, \rbnum_{\infty,\infty}
(\pdeg,\xarity)$,
\end{enumerate}
where $\rbnum_{1, \infty}(\pdeg, \xarity)$, $\rbnum_{\infty,
1}(\pdeg, \xarity)$, and  $\rbnum_{\infty, \infty}(\pdeg, \xarity)$
are given by {\rm Corollary~\mref{thm:five}}.
\end{theorem}
\nc\Integer{\NN}
\proofbegin We will prove only (3).
Let $w_\infty \in \rbset_{\infty, \infty}(\pdeg, \xarity;
\prun, \xrun)$.
Each element $w \in \rbset_{\inftyv,\inftyv}(\pdeg,\xarity; \prun,
\xrun)$ whose monochromatic image is $w_\infty$ is obtained from
$w_\infty$ by coloring the $P$-\runs and
$x$-\runs in $w_\infty$. For any positive integer $j$, let
$\Integer_j$ be the set consisting of the first $j$ natural numbers.
Since the number of mappings from $\Integer_\pdeg$
to $\Integer_\pnumb$ is $\pnumb^\pdeg$ and the number of mappings from
$\Integer_\xarity$ to $\Integer_\xnumb$ is $\xnumb^\xarity$, the
result follows.
Alternatively, this theorem also follows by summing over all $\prun$
and $\xrun$ the cardinalities from Eq.~(\mref{eq:colcol}) and using
Corollary \mref{thm:five}.  \proofend

\begin{coro} \mlabel{thm:corol11}
The generating functions for the double sequences $\rbnum_{\vec{1},
\inftyv}(\pdeg, \xarity)$, $\rbnum_{\inftyv, \vec{1}}(\pdeg,
\xarity)$, and $\rbnum_{\inftyv, \inftyv}(\pdeg,
\xarity)$ are, respectively,
\renewcommand{\theenumi}{\arabic{enumi}}  
\renewcommand{\theenumii}{\alph{enumii}}  
\renewcommand{\labelenumi}{{\rm (\arabic{enumi})}}  
\renewcommand{\labelenumii}{{\rm (\alph{enumii})}}  
\begin{enumerate}
\item $\displaystyle \rbfun_{\vec{1}, \inftyv}(\pdegvar, \xarityvar) =
\rbfun_{1,1}\left(\pdegvar, \, \frac{\xnumb \xarityvar}{1 - \xnumb
\xarityvar}\right)$, \vspace{0.1in}
\item $\displaystyle
\rbfun_{\inftyv, \vec{1}}(\pdegvar, \xarityvar) =
\rbfun_{1,1}\left(\frac{\pnumb \pdegvar}{1 -
\pnumb \pdegvar}, \xarityvar \right)$, and \vspace{0.1in}
\item $\displaystyle  \rbfun_{\inftyv, \inftyv}(\pdegvar,
\xarityvar) = \rbfun_{1,1}\left(\frac{\pnumb \pdegvar}{1 -
\pnumb \pdegvar}, \frac{\xnumb \xarityvar}{1 - \xnumb
\xarityvar}\right)$.
\end{enumerate}
\end{coro}

\proofbegin
By the above theorem and Corollary \mref{thm:seven}, for case (3), we
have
\begin{equation*}
\begin{array}{rcl}
\rbfun_{\inftyv, \inftyv}(\pdegvar, \xarityvar)  &=&
\displaystyle \sum_{\pdeg=1}^\infty
\sum_{\xarity=1}^\infty \pnumb^\pdeg \xnumb^\xarity \rbnum_{\infty,
\infty}(\pdeg, \xarity)\pdegvar^\pdeg \xarityvar^\xarity
\vspace{0.2in}\\ &=& \displaystyle \sum_{\pdeg=1}^\infty
\sum_{\xarity=1}^\infty \rbnum_{\infty,
\infty}(\pdeg, \xarity)(\pnumb\pdegvar)^\pdeg
(\xnumb \xarityvar)^\xarity \vspace{0.2in}\\
&=&\displaystyle
\rbfun_{1,1}\left(\frac{\pnumb \pdegvar}{1 -
\pnumb \pdegvar}, \frac{\xnumb \xarityvar}{1 - \xnumb
\xarityvar}\right). \ \square\end{array}
\end{equation*}

\subsection{General case}
\mlabel{sec:threeptthree}

\begin{theorem} \mlabel{thm:ten}
Let $\rbset_{\pexpv, \xexpv}$ be the set of Rota-Baxter words
with $\pnumb$ operators $\pv$ having \exponent vector $\pexpv$, and
$\xnumb$ generators $\xv$ having \exponent vector $\xexpv$, where
$\pexpv$ \paren{resp.~$\xexpv$} may have finite or infinite
components. Let $\colorGF_{\pnumb,\pexpv}(\pdegvar)$
\paren{resp.~$\colorGF_{\xnumb,\xexpv}(\xarityvar)$}
be the
generating
function for the sequence $\colornum(\pa,\pnumb,\pexpv)$ as $\pa$
varies \paren{resp.~$\colornum(\xb,\xnumb,\xexpv)$ as $\xb$ varies}.
 Then the
generating function for $\rbnum_{\pexpv, \xexpv}(\pdeg, \xarity)$ is
$$\rbfun_{\pexpv,\xexpv}(\pdegvar,\xarityvar)
=\rbfun_{1,1}(\colorGF_{\pnumb,\pexpv}(\pdegvar),
\colorGF_{\xnumb,\xexpv}(\xarityvar)).
$$
\end{theorem}

\proofbegin The special cases when either $\pnumb=1$, or $\xnumb=1$,
or both $\pnumb=\xnumb=1$ specialize to corresponding cases of Theorem
\mref{thm:six}.  Thus we shall assume $\pnumb \geq 2$ and $\xnumb \geq
2$.  Some special cases when $\pexpv = \onev$ or $\inftyv$ and when
$\xexpv = \onev$ or $\inftyv$ are included in Corollaries
\mref{thm:seven} and \mref{thm:corol11}.

Consider the
forgetful map
$\pxmapscalar$ of Eq.~(\mref{eq:pxmapscalar}), which is surjective:
\begin{equation*}
\pxmapscalar:  \rbset_{\pexpv, \xexpv}(\pdeg, \xarity; \prun, \xrun)
\longrightarrow \rbset_{\infty, \infty}(\pdeg, \xarity; \prun, \xrun)
\end{equation*}

We have obviously
$\rbset_{\pexpv, \xexpv}(\pdeg, \xarity; \prun, \xrun) =
(\pxmapscalar)^{-1}\left(\rbset_{\infty, \infty}(\pdeg, \xarity;
\prun, \xrun)\right)$, and using the bijection from
Eq.~(\mref{eq:colcol}), we have
\begin{equation*}
\begin{array}{rcl}
&&r_{\pexpv,\xexpv}(\pdeg,\xarity; \prun, \xrun)\\&=&\displaystyle
 r_{1,1}(\prun,\xrun)
\left(\sum_{\pdegv\, \in \comset(\pdeg, \prun, \infty)}
\prod_{i=1}^\prun \colornum(\pdeg_i, \pnumb, \pexpv)
\right)
\left(\sum_{\xarityv \,\in \comset(\xarity, \xrun, \infty)}
\prod_{j=1}^\xrun \colornum(\xarity_j, \xnumb, \xexpv)
\right)
\end{array}
\end{equation*}

Recalling Eq.~(\mref{eq:rundecomp}), we now calculate the generating
function.
\begin{equation*}
\begin{array}{rcl}
&&\displaystyle
\rbfun_{\pexpv,\xexpv}(z,t)= \sum_{\pdeg\geq 1} \sum_{\xarity\geq 1}
r_{\pexpv,\xexpv}(\pdeg,\xarity) z^\pdeg t^\xarity
\vspace{0.15in}\\
&=&\displaystyle
\sum_{\pdeg, \xarity, \prun, \xrun}
r_{1,1}(\prun,\xrun)
\left(\sum_{\pdegv\, \in \comset(\pdeg, \prun, \infty)}
\prod_{i=1}^\prun \colornum(\pdeg_i, \pnumb, \pexpv)
\pdegvar^\pdeg
\right)
\left(\sum_{\xarityv \,\in \comset(\xarity, \xrun, \infty)}
\prod_{j=1}^\xrun \colornum(\xarity_j, \xnumb, \xexpv)
\xarityvar^\xarity
\right) \vspace{0.15in}\\
&=&\displaystyle
\sum_{\prun\geq 1}\sum_{\xrun\geq 1}
r_{1,1}(\prun,\xrun)
\left(\sum_{\pdeg\geq 1}
\sum_{\pdegv} 
\prod_{i=1}^\prun \colornum(\pdeg_i, \pnumb, \pexpv)
\pdegvar^\pdeg
\right)
\left(\sum_{\xarity\geq 1}
\sum_{\xarityv} 
\prod_{j=1}^\xrun \colornum(\xarity_j, \xnumb, \xexpv)
\xarityvar^\xarity
\right)\,.\\
\end{array}
\end{equation*}

Working
on the last parenthesized expression,
we have \begin{equation*}
\begin{array}{rcl}
\displaystyle
\sum_{\xarity\geq 1} \left[\sum_{\xarityv \,\in
\comset(\xarity,
\xrun, \infty)} \prod_{j=1}^\xrun \colornum(\xarity_j, \xnumb, \xexpv)
\right] \xarityvar^\xarity
&=&\displaystyle
\sum_{\xarity\geq 1} \left[\sum_{\xarityv \,\in
\comset(\xarity,
\xrun, \infty)} \prod_{j=1}^\xrun \colornum(\xarity_j, \xnumb, \xexpv)
\xarityvar^{\xarity_j} \right]
\vspace{0.15in}\\
&=&\displaystyle
\left (\sum_{\xb\geq 1} \colornum(\xb, \xnumb, \xexpv)
\xarityvar^{\xb} \right )^\xrun
\vspace{0.15in}\\
&=&\displaystyle
\bigl(\colorGF_{\xnumb, \xexpv}(\xarityvar)\bigr)^\xrun\,.
\end{array}
\end{equation*}
The result now follows from the definition of $\rbfun_{1,1}(\pdegvar,
\xarityvar)$. \ $\square$

\begin{coro} If the \exponent vectors $\pexpv, \xexpv$ are uniform,
say $\pexpv = \vec{\prunlen}$ and $\xexpv = \vec{\xrunlen}$ for some
finite or infinite constants $\prunlen, \xrunlen$, then the generating
function for $\rbnum_{\vec{\prunlen},\vec{\xrunlen}}(\pdeg, \xarity)$
is:
$$\rbfun_{\vec{\prunlen},\vec{\xrunlen}}(\pdegvar,\xarityvar)
=\rbfun_{1,1}\left ( \frac{\pnumb \, \comGF_{1, \prunlen}( \pdegvar)}
        {1-(\pnumb-1)\comGF_{1, \prunlen}( \pdegvar)}\,,\,
    \frac{\xnumb\, \comGF_{1, \xrunlen}(\xarityvar)}
        {1-(\xnumb-1)\comGF_{1, \xrunlen}(\xarityvar)}\right ).$$
\end{coro}

\proofbegin
This follows from
Lemma \mref{thm:lemma8} and Corollary \mref{thm:corol9}.\ $\square$
\deleted{
 since
$$\colorGF_{\xnumb, \xexpv}(\xarityvar)
=
 \frac{\xnumb \comGF_{1, \xrunlen}(\xarityvar)}
{1 - (\xnumb-1) \comGF_{1, \xrunlen}(\xarityvar)}\,.\ \square
$$
}

\vspace{-0.1in}
\section{Conclusion and outlook} \mlabel{sec:four}

\vspace{-0.1in}
We have obtained generating functions and algorithms related to
enumeration of sets of Rota-Baxter words
in various generalities.
It is interesting to see the close relation with Catalan numbers,
further revealing the combinatorial nature  of Rota-Baxter algebras in
cases where the sets form canonical bases.

The cases we have considered allow multiple operators and multiple
generators with uniform exponents. The case where the operators and
generators are allowed to have variable exponents are in principle
solved and it only remains to compute the generating functions for the
number of colorings. Since the coloring problem is a
constraint satisfaction problem, we are investigating this
connection to see if the result is alrady known.   We expect these
results to be
useful in the construction of more general Rota-Baxter type algebras,
where the unary operators may satisfy other identities instead.

We mentioned a few new integer sequences in this study.  Clearly,
since the most general generating functions given by Theorem
\mref{thm:ten} are parametrized by two arbitrary vectors, we expect
many specialized sequences to be unknown in the Sloane database.  The
connections between sets of Rota-Baxter words and other sets of
combinatorial objects when their counting sequences are the same will
also be of tremendous interest in expanding the realm of symbolic
computation because such connections will allow the representation of
these other sets of combinatorial objects using Rota-Baxter words
which are purely algebraic in nature.

\medskip

\noindent
{\bf Acknowledgements:  } The authors acknowledge support from NSF
grants DMS 0505643 (Li Guo) and NSF grants CCF-0430722 (William
Y.~Sit), and thank K.~Ebrahimi-Fard and W.~Moreira for helpful
discussions.

\vfill

%
%
\bibliographystyle{elsart-harv}


\end{document}